\documentclass[12pt,leqno]{article}
\usepackage{amsfonts}
\usepackage{amsmath, amsthm, amsfonts, amssymb, color}
\usepackage{mathrsfs}
\usepackage{color}
\theoremstyle{definition}

\newcommand{\scr}[1]{\mathscr #1}
\definecolor{wco}{rgb}{0.5,0.2,0.3}

\numberwithin{equation}{section} \theoremstyle{remark}

\newcommand{\ua}{\uparrow}

\title{{\bf   Singular Density Dependent    Stochastic Differential Equations }\footnote{Supported in
 part by  the National Key R\&D Program of China (No. 2022YFA1006000, 2020YFA0712900) and NNSFC (11831014, 11921001).} }
\author{{\bf  Feng-Yu Wang   }\\
\footnotesize{ Center for Applied Mathematics, Tianjin University, Tianjin 300072, China}\\
  \footnotesize{  Department of Mathematics,
Swansea University, Bay Campus, SA1 8EN, United Kingdom}\\
\footnotesize{    wangfy@tju.edu.cn}}
\begin{document}
\allowdisplaybreaks
\def\R{\mathbb R}  \def\ff{\frac} \def\ss{\sqrt} \def\B{\mathbf
B}
\def\N{\mathbb N} \def\kk{\kappa} \def\m{{\bf m}}
\def\ee{\varepsilon}\def\ddd{D^*}
\def\dd{\delta} \def\DD{\Delta} \def\vv{\varepsilon} \def\rr{\rho}
\def\<{\langle} \def\>{\rangle}
  \def\nn{\nabla} \def\pp{\partial} \def\E{\mathbb E}
\def\d{\text{\rm{d}}} \def\bb{\beta} \def\aa{\alpha} \def\D{\scr D}
  \def\si{\sigma} \def\ess{\text{\rm{ess}}}\def\s{{\bf s}}
\def\beg{\begin} \def\beq{\begin{equation}}  \def\F{\scr F}
\def\Ric{\mathcal Ric} \def\Hess{\text{\rm{Hess}}}
\def\e{\text{\rm{e}}} \def\ua{\underline a} \def\OO{\Omega}  \def\oo{\omega}
 \def\tt{\tilde}\def\[{\lfloor} \def\]{\rfloor}
\def\cut{\text{\rm{cut}}} \def\P{\mathbb P} \def\ifn{I_n(f^{\bigotimes n})}
\def\C{\scr C}      \def\aaa{\mathbf{r}}     \def\r{r}
\def\gap{\text{\rm{gap}}} \def\prr{\pi_{{\bf m},\nu}}  \def\r{\mathbf r}
\def\Z{\mathbb Z} \def\vrr{\nu} \def\ll{\lambda}
\def\L{\scr L}\def\Tt{\tt} \def\TT{\tt}\def\II{\mathbb I}
\def\i{{\rm in}}\def\Sect{{\rm Sect}}  \def\H{\mathbb H}
\def\M{\mathbb M}\def\Q{\mathbb Q} \def\texto{\text{o}} \def\LL{\Lambda}
\def\Rank{{\rm Rank}} \def\B{\scr B} \def\i{{\rm i}} \def\HR{\hat{\R}^d}
\def\to{\rightarrow} \def\gg{\gamma}
\def\EE{\scr E} \def\W{\mathbb W}
\def\A{\scr A} \def\Lip{{\rm Lip}}\def\S{\mathbb S}
\def\BB{\scr B}\def\Ent{{\rm Ent}} \def\i{{\rm i}}\def\itparallel{{\it\parallel}}
\def\g{{\mathbf g}}\def\Sect{{\mathcal Sec}}\def\T{\mathcal T}\def\BB{{\bf B}}
\def\f\ell \def\g{\mathbf g}\def\BL{{\bf L}}  \def\BG{{\mathbb G}}
\def\Bd{{D^E}} \def\BdP{D^E_\phi} \def\Bdd{{\bf \dd}} \def\Bs{{\bf s}} \def\GA{\scr A}
\def\Bg{{\bf g}}  \def\Bdd{\psi_B} \def\supp{{\rm supp}}\def\div{{\rm div}}
\def\ddiv{{\rm div}}\def\osc{{\bf osc}}\def\1{{\bf 1}}\def\BD{\mathbb D}
\def\H{{\bf H}}\def\gg{\gamma} \def\n{{\mathbf n}}
\maketitle

\begin{abstract} The (strong and weak) well-posedness is proved for singular SDEs depending on the distribution density  point-wisely and globally, where the drift satisfies a
local integrability condition in time-spatial variables, and is Lipschitz continuous in the distribution density with respect to a local $L^k$-norm. Density dependent reflecting SDEs are also studied.
  \end{abstract} \noindent

 \vskip 2cm

 \section{Introduction}
 The study of distribution dependent SDEs goes back to McKean's pioneering work \cite{MC} where an expectation dependent SDE is proposed to characterize Maxwellian gas.
 Comparing with the dependence on the global distribution, the point-wise dependence on the density function is more singular for SDEs.
The  point-wisely  density dependent SDE is called Nemytskii-type McKean-Vlasov SDE,  see \cite{BR1, BR2} for the correspondence of this type SDEs and nonlinear PDEs.

In recent years, distribution dependent SDEs have  been intensively investigated and a plenty of results have been derived.    However, much less is known for density dependent SDEs, see Remark 1.1 below for existing results.

Let $\ell_\xi: \R^d\to [0,\infty)$   be  the distribution density function of an absolutely continuous random variable $\xi$ on $\R^d$.  We investigate  the following   SDE    depending on the distribution density   point-wisely and globally:
 \beq\label{E1} \d X_t= b_t(X_t,\ell_{X_t}(X_t), \ell_{X_t})\d t +\si_t(X_t,\ell_{X_t})\d W_t,\ \ \ \ t\in [0,T],\end{equation}
where  $T>0$ is fixed,   $\{W_t\}_{t\in [0,T]}$ is an $m$-dimensional Brownian motion on a complete filtration probability space $(\OO, \{\F_t\}_{t\in [0,T]}, \P)$, and
$$b: [0,T]\times \R^d\times [0,\infty)\times \D_+^1\to \R^d,\ \ \si: [0,T]\times\R^d\times \D_+^1\to \R^d\otimes\R^m$$ are measurable, where
  $$\D_+^1:= \bigg\{f\in L^1(\R^d): f\ge 0,   \int_{\R^d} f(x)\d x\le 1 \bigg\} $$ is a closed subspace of  $L^1(\R^d)$.

\beg{defn}  A continuous adapted process $(X_t)_{t\in [0,T]}$ on $\R^d$  is called a (strong) solution of \eqref{E1}, if
$$ \int_0^T \E\big[|b_s(X_s,\ell_{X_s}(X_s), \ell_{X_s})|+\|\si_s(X_s,  \ell_{X_s})\|^2\big]  \d s<\infty $$ and
$\P$-a.s.
$$X_t=X_0+ \int_0^t b_s(X_s,\ell_{X_s}(X_s), \ell_{X_s})\d s +\int_0^t \si_s(X_s,  \ell_{X_s})\d W_s, \ \  t\in [0,T].$$

  A pair $(X_t,W_t)_{t\in [0,T]}$ is called a weak solution of \eqref{E1}, if $(W_t)_{t\in [0,T]}$ is an $m$-dimensional Brownian under a complete filtration probability space
$(\OO, \{\F_t\}_{t\in [0,T]}, \P)$ such that $(X_t)_{t\in [0,T]}$ solves $\eqref{E1}$. We identify   any two weak solutions $(X_t, W_t)_{t\in [0,T]}$ and $(\bar X_t,\bar W_t)_{t\in [0,T]}$ if
$(X_t)_{t\in [0,T]}$ and $(\bar X_t)_{t\in [0,T]}$  have the same   distribution under  the corresponding probability spaces.  \end{defn}

\paragraph{Remark 1.1.} We introduce below some existing results concerning the well-posdedness of \eqref{E1} in two special situations.

 (1) When $m=d,\si= I_d$  (the $d\times d$ identity matrix), and $b_t(x,r,\rr)= b_t(x,r)$ does not depend on $\rr$, the weak solutions are studied in \cite{HRZ, IR}. In  \cite{HRZ},
  the weak existence is proved for  $b_t(x,r)$  bounded and continuous in $(t,r)$ locally uniformly in $x$,  and the weak and strong uniqueness holds when   $b_t(x,r)$ is furthermore Lipschitz continuous in $r$ uniformly in $(t,x)$.  In \cite{IR} the initial density is in $C^{\bb+}:=\cup_{p>\bb} C^p$ for some $\bb\in (0,\ff 1 2)$,   the weak well-posedness is proved for  $b_t(x,r):= F(r) \tt b_t(x)$,  where  $\tt b\in C([0,T];C^{-\bb})$ and  $F$  is bounded and Lipschitz continuous such that $rF(r)$ is Lipschitz continuous in $r\ge 0$. See \cite{RU1} and references within for the case with better  drift.

(2) In \eqref{E1} the noise does not point-wisely depend on the density. It seems that to solve  SDEs with point-wisely density dependent noise, one needs stronger regularity    for   the initial density and  the coefficients.  For instance, \cite{RU2} proved the well-posedness and studied  the  propagation of chaos for the following  SDE with point-wisely density dependent noise:
$$\d X_t= b(\ell_{X_t}(X_t))\d t+ \si(\ell_{X_t}(X_t)) \d W_t,$$
where the initial distribution density is   $C^{2+}$-smooth, $b$ is   $C^2$-smooth, and $\si$ is uniformly elliptic and $C^3$-smooth.

In this paper, we only consider  density dependent
SDEs with singular (non-continuous   in spatial) drifts, and leave  to a forthcoming paper for  the study of  regular SDEs with point-wisely density dependent  noise.

   \

 In Section 2, we state two main results of the paper which provide the well-posedness of \eqref{E1} for     density free noise and density dependent noise respectively,  and explain the main idea of the proof. To realize the idea, in Section 3 we recall some heat kernel estimates based on \cite{MPZ} and present new estimates.  With these estimates we  prove the main results in  Sections 4 and 5 respectively, and finally  make an extension to the reflecting setting in Section 6.

 \section{Main results and idea of proof}

 To characterize the time-spatial singularity, we   recall some spaces of locally integrable functions  introduced in \cite{XXZZ}.

For $p \in [1,\infty]$ and    $f\in\B(\R^d)$, the space of measurable functions on $\R^d$, let
$$\|f\|_{L^p}:= \bigg(\int_{\R^d} |f|^p(x)\d x\bigg)^{\ff 1 p},\ \ \|f\|_{\tt L^p} := \sup_{z\in \R^d}  \|1_{B(z,1)}f\|_{L^p}<\infty,$$
where $B(z,r):= \{x\in\R^d: |x-z|\le r\}, r>0, z\in \R^d.$ We write $f\in L^p\ (f\in \tt L^p)$ if $\|f\|_{L^p}<\infty\ (\|f\|_{\tt L^p}<\infty)$.

For any $p,q\in [1,\infty]$ and   $f\in \B([0,T]\times\R^d)$, the space of measurable functions on $[0,T]\times \R^d$, let
$$\|f\|_{L_q^p} := \bigg(\int_0^T \|f_t\|_{L^p}^q\d t\bigg)^{\ff 1 q},\ \ \|f\|_{\tt L_q^p} := \sup_{z\in \R^d} \bigg(\int_0^T \|1_{B(z,1)} f_t\|_{L^p}^q\d t\bigg)^{\ff 1 q}.$$
We denote $f\in L_q^p (f\in \tt L_q^p)$ if $\|f\|_{L_q^p}<\infty\ (\|f\|_{\tt L_q^p}<\infty)$.

In the following  the parameter $(p,q)$ will be taken from the class
$$\scr K:=\Big\{(p,q)\in (2,\infty]: \ff d p+\ff 2 {q}<1\Big\}.$$

For simplicity, we identify $L^\infty =\tt L^\infty$ with $\B_b(\R^d)$, the space of bounded measurable functions on $\R^d$, equipped with the uniform norm
  $$ \|f\|_{L^\infty}=\|f\|_{\tt L^\infty}=\|f\|_\infty:=\sup_{\R^d} |f|.$$   This uniform norm is defined for real functions on an abstract space.
   Similarly, $\tt L_\infty^\infty= L_\infty^\infty=\B_b([0,T]\times\R^d)$ is the space   of bounded measurable functions on $[0,T]\times \R^d$ equipped with the uniform norm,
and  $L_\infty^k\ (\tt L_\infty^k)$ is the space of  functions $f\in\B([0,T]\times\R^d)$ such that
  $$\|f\|_{L_\infty^k}:=\sup_{t\in [0,T]} \|f_t\|_{L^k}<\infty\ \big(\|f\|_{\tt L_\infty^k}:=\sup_{t\in [0,T]} \|f_t\|_{\tt L^k}<\infty\big).$$
 Finally,  let $\nn$ be the gradient in $\R^d$, and let $\|\nn f\|_\infty$ denote the Lipschitz constant of a real function $f$ on $\R^d$.

In the following, we state our main results for density free noise and density   dependent noise respectively, and briefly explain the main idea of proof.

 \subsection{Density free noise}

 In this part, we let $\si_t(x,\rr)=\si_t(x)$ do  not depend on $\rr$. For $k>1$ and a signed measure $\mu$ with density function $\ell_\mu(x):=\ff{\mu(\d x)}{\d x}$, let
 $$\|\mu\|_{L^k}:=  \| \ell_\mu\|_{L^k},\ \ \|\mu\|_{\tt L^k}:=  \|\ell_\mu \|_{\tt L^k}.$$
 When $k=1$,   we define
 $$\|\mu\|_{L^1}:=\sup_{\|f\|_\infty\le 1} |\mu(f)|,\ \ \   \|\mu\|_{\tt L^1}:=\sup_{z\in\R^d} \sup_{\|f\|_\infty\le 1} |\mu(1_{B(z,1)}f)|,$$ where $\mu(f):=\int_{\R^d}f\d\mu.$
   Note that $\|\cdot\|_{L^1}$ is the total variation norm.

Let $\scr P$ be the set of all probability measures on $\R^d$. We will solve  \eqref{E1} with
 initial distributions in the classes
 $$  {\scr P}^k:=   \Big\{\nu\in \scr P: \|\nu\|_{L^k} <\infty\Big\}, \ \ \
    \tt{\scr P}^k:=\Big\{\nu\in \scr P: \|\nu\|_{\tt L^k} <\infty\Big\},\ \ k\in [1,\infty],   $$
  which are complete metric spaces under distances
    $\|\nu_1-\nu_2\|_{L^k}$ and $\|\nu_1-\nu_2\|_{\tt L^k}$ respectively.

 \beg{enumerate}  \item[{\bf (A)}] $a_t(x):= (\si_t\si_t^*)(x)$ and $b_t(x,r,\rr)= b_t^{(1)}(x)+ b_t^{(0)}(x,r,\rr)$ satisfy    the following conditions
 for some $k\in [1,\infty]$.
 \item [$(A_1)$] $a_t(x) $ is invertible with   $\|a\|_\infty+\|a^{-1}\|_\infty<\infty$,  and there exist constants  $\aa\in (0,1)$ and $C>0$ such that
 $$\sup_{t\in [0,T]} \|a_t(x)-a_t(y)\|\le C|x-y|^\aa,\ \ x,y\in \R^d.$$
 \item[$(A_2)$]  There exist  $(p_0,q_0)\in \scr K$,  $\theta>   \ff{2}{q_0} +\ff d{p_0} -1,$ and $1\le f_0\in \tt L_{q_0}^{p_0}$ such that
  \beg{align*} &  |b_t^{(0)}(x,r,\rr)- b_t^{(0)}(x,\tt r,\tt\rr)|\le  f_0(t,x) t^\theta \big(|r-\tt r|+\|\rr-\tt\rr\|_{\tt L^k}\big),\\
 & |b_t^{(0)}(x,r,\rr)|\le f_0(t,x),\ \     \   (t,x)\in (0,T]\times\R^d,  r,\tt r\in [0,\infty), \rr,\tt\rr\in \tt L^k\cap \D_+^1.\end{align*}
 \item[$(A_3)$] $b^{(1)}_t(0)$  is bounded in $t\in [0,T]$ and
 $$\|\nn b^{(1)}\|_\infty:=\sup_{t\in [0,T]} \sup_{x\ne y} \ff{|b_t^{(1)}(x)-b_t^{(1)}(y)|}{|x-y|}<\infty.$$
 \end{enumerate}

To ensure   $\ell_{X_\cdot}\in L^k$ for $\ell_{X_0}\in L^k$, we replace  $(A_2)$ by the following $(A_2')$.
 \beg{enumerate}   \item[$(A_2')$]  There exist  $C\in (0,\infty),  (p_0,q_0)\in \scr K$,  $\theta>   \ff{2}{q_0} +\ff d{p_0} -1,$     and $0\le f_0\in   L_{q_0}^{p_0}$ such that
  \beg{align*} &  |b_t^{(0)}(x,r,\rr)- b_t^{(0)}(x,s,\tt\rr)|\le  t^\theta \big(C+f_0(t,x)\big)\big(|r-s|+\|\rr-\tt\rr\|_{L^k}\big),\\
 & |b_t^{(0)}(x,r,\rr)|\le f_0(t,x),\ \     \   (t,x)\in (0,T]\times\R^d,  r,s\in [0,\infty), \rr,\tt\rr\in   L^k\cap \D_+^1.\end{align*} \end{enumerate}

Under the above assumptions, the following result ensures the well-posedness of \eqref{E1} for initial distributions in
$\tt{\scr P}^k$ or $\scr P^k$ for $$k\in \Big[\ff{p_0}{p_0-1},\infty\Big]\cap (k_0,\infty],\ \
 k_0:=\ff d{2\theta+1-2q_0^{-1} - dp_0^{-1}}.$$
This explains the role played by the quantity  $\theta$  in $(A_2)$ and $(A_2')$: for bigger $\theta$,
  $(A_2)$ and $(A_2')$ provide stronger upper bound condition on  $|b_t^{(0)}(x,r,\rr)-b_t^{(0)} (x,\tt r,\tt\rr)|$   for small $t$, so that
the SDE is solvable  for initial distributions in larger classes   $\scr P^k$ and $\tt{\scr P}^k$.
 In particular, when $p_0=\infty$ and $\theta$ is large enough such that $k_0<1$, we may take $k=1$ so that the SDE is well-posed for any initial distribution $\nu\in \scr P.$

  \beg{thm}\label{T1} Let    $k\in [\ff{p_0}{p_0-1},\infty]$ with  $k>k_0:= \ff d{2\theta+1-2q_0^{-1} - dp_0^{-1}}.$
\beg{enumerate} \item[$(1)$]  Under {\bf (A)},  for any   $\nu\in \tt {\scr P}^k,$ $\eqref{E1}$ has a unique  weak  solution with $\L_{X_0}=\nu$ satisfying $\ell_{X_\cdot}\in \tt L_\infty^k$,  and  there exist an increasing function $\LL: [0,\infty)\to (0,\infty)$   such that for  any two  weak  solutions  $\{X_t^i\}_{i=1,2}$  of $\eqref{E1}$ with
  $\ell_{X_\cdot^i} \in\tt L^k_\infty$,
 \beq\label{S2} \sup_{t\in [0,T]} \|\ell_{X_t^1} -\ell_{X_t^2}\|_{\tt L^k}\le \LL \big(\|\L_{X_0^1}\|_{\tt L^k}\land  \|\L_{X_0^1}\|_{\tt L^k}\big) \|\L_{X_0^1}- \L_{X_0^2}\|_{\tt L^k}.\end{equation}
 If  moreover $\si_t$ is weakly differentiable with
 \beq\label{SI}  \|\nn\si\|\le \sum_{i=1}^l f_i\ \text{for\  some\ }
  l\in \mathbb N, 0\le f_i\in \tt L_{q_i}^{p_i}, (p_i,q_i)\in \scr K, 1\le i\le l,\end{equation}
  then for any   $X_0$ with   $\L_{X_0}\in  \tt L^k$, $\eqref{E1}$ has a unique  strong  solution with $\ell_{X_\cdot}\in \tt L_\infty^k.$
 \item[$(2)$] Under  {\bf (A)} with $(A_2')$ replacing $(A_2)$,
   assertions in   $(1)$  hold  for $(\scr P^k, L_\infty^k, L^k)$ replacing $(\tt{\scr P}^k, \tt L_\infty^k, \tt L^k)$.\end{enumerate}
  \end{thm}

   \subsection{Density dependent noise}

 In this part, we allow $\si$ to be density dependent but make stronger assumptions on the initial density and the coefficients in the spatial variable.

 For any $n\in \mathbb Z^+$, let
 $C_b^n(\R^d)$ be the class of real functions $f$ on $\R^d$ with continuous derivatives $\{\nn^i f\}_{0\le i\le n}$ such that
 $$\|f\|_{C_b^n}:= \sum_{i=0}^n \|\nn^i f\|_\infty<\infty.$$
 For any $n\in \mathbb Z^+$ and $\aa\in (0,1)$, $C_b^{n+\aa}(\R^d)$ is the space of functions $f\in C_b^n(\R^d)$ such that
 $$\|f\|_{C_b^{n+\aa}}:= \|f\|_{C_b^n} + \sup_{x\ne y} \ff{|\nn^n f(x)-\nn^n f(y)|}{|x-y|^\aa}<\infty.$$

 \beg{enumerate} \item[{\bf(B)}]    There exist $1\le f_0\in\tt L_{q_0}^{p_0}$,   $ C\in (0,\infty) $ and $\aa\in (0,1)$, such that
 the following conditions hold for all $t\in (0,T],\ x,y\in \R^d, \ r,\tt r\in [0,\infty)$ and $  \rr,\tt\rr\in  L^\infty\cap \D_+^1$:
$$|b_t(x,r,\rr)|\le f_0(t,x),$$
$$    |b_t(x,r,\rr)-b_t(x,\tt r,\tt\rr)|\le C  (|r-\tt r|+  \|\rr-\tt\rr\|_{\infty}),$$
  $$  \|\si\|_\infty+\|\nn\si\|_\infty+\|(\si\si^*)^{-1}\|_\infty\le C,$$
$$    \|\nn \si_t(\cdot, \rr)(x)-\nn\si_t(\cdot,  \rr)(y)\| \le C |x-y|^\aa,$$
$$      \| \si_t(\cdot,\rr)- \si_t(\cdot,\tt\rr) \|_{C_b^\aa}
     \le   C \|\rr-\tt\rr\|_{\infty}.$$
   \end{enumerate}

  \beg{thm}\label{T2} Assume {\bf (B)} and let $\bb\in (0,1-\ff{d}{p_0}-\ff 2 {q_0})$.   For any initial value $($initial density$)$ with $\ell_{X_0}\in C_b^\bb(\R^d)$,
   $\eqref{E1}$ has a unique  strong  $($weak$)$ solution satisfying $\ell_{X_\cdot}\in L_\infty^\infty$, and there exists a constant $c>0$ such that
   \beq\label{A1} \sup_{t\in [0,T]}\|\ell_{X_t}\|_{C_b^\bb} \le  c \|\ell_{X_0}\|_{C_b^\bb}. \end{equation}
Moreover, there exists   an increasing function $\LL: (0,\infty)\to (0,\infty)$ such that   for any two  solutions $\{X_t^i\}_{i=1,2}$ with $\ell_{X_0^i}\in C_b^\bb(\R^d)$ and $\ell_{X_\cdot^i}\in  L_\infty^\infty$,
    \beq\label{A2} \sup_{t\in [0,T]}\|\ell_{X_t^1}-\ell_{X_t^2}\|_{\infty} \le    \LL\big(\|\ell_{X_0^1}\|_{C_b^\bb}\land \|\ell_{X_0^1}\|_{C_b^\bb}\big)\|\ell_{X_0^1}-\ell_{X_0^2}\|_\infty.
     \end{equation}
    \end{thm}

 \subsection{Idea of proof}
 For fixed  $k\ge 1$  and   $\nu\in \tt{\scr P}^k,$ let  $\tt{\scr P}^k_{\nu,T}$   be the set of all bounded measurable maps
 $$\gg: (0,T]\to \tt L^k\cap \D_+^1,\ \ \gg_0=\nu.$$ When $k=1,$ the initial value
   $\gg_0$ may be singular, and if it is absolutely continuous we regard it as its density function.

 Then $\tt{\scr P}^k_{\nu,T}$ is complete under the metric
 $$\tt d_{k,\ll}(\gg^1,\gg^2):= \sup_{t\in [0,T]}\e^{-\ll t} \|\gg^1_t-\gg_t^2\|_{\tt L^k},\ \ \gg^1,\gg^2\in
 \tt{\scr P}^k_{\nu,T}$$  for $\ll>0$.
We define  $({\scr P}^k_{\nu,T}, d_{k,\ll})$ in the same way with $(L^k, \scr P^k)$ replacing
$(\tt L^k, \tt{\scr P}^k).$

  For any $\gamma \in \tt{\scr P}^k_{\nu,T}$, let
\beg{align*}b_t^\gg(x):= b_t(x,\gg_t(x),\gg_t),\ \ \si_t^\gg(x):= \si_t(x,  \gg_t),\ \ t\in (0,T], x\in\R^d.\end{align*} Then  for $\nu:=\L_{X_0}\in \tt L^k$,  \eqref{E1} has a unique   (weak or strong) solution with $\ell_{X_\cdot}\in \tt L^k_\infty$
  if we could verify
  the following two things:
\beg{enumerate}\item[$1)$] For any $\gg\in \tt{\scr P}_{\nu,T}^k,$ the SDE
\beq\label{E2} \d  X_t^\gg = b_t^\gg(X_t^\gg)\d t+\si_t^\gg(X_t^\gg)\d W_t,\ \ t\in [0,T], \ X_0^\gg=X_0 \end{equation}
is (weakly or strongly) well-posed, and
$$\gg\mapsto \Phi^{\nu}_t\gg:=  \ell_{X_t^\gg},\ \ t\in (0,T],\ \Phi^{\nu}_0\gg:=\gg_0=\nu$$
provides a map   $\Phi^{\nu}: \tt{\scr P}^k_{\nu,T} \to \tt{\scr P}^k_{\nu,T}.$
\item[$2)$] $\Phi^{\nu}$ has a unique fixed point $\bar\gg$ in $\tt{\scr P}^k_{\nu,T}$.  \end{enumerate}
Indeed, from  these  we see that  $X_t:=X_t^{\bar \gg}$ is the unique (weak or strong) solution of \eqref{E1}
with $\L_{X_\cdot}\in \tt L_\infty^k$.

\

To verify 1) and 2), in Section 2 we  recall some heat kernel upper bounds  of \cite{MPZ},
and estimate the $\tt L^p_q$-$\tt L^{p'}_q$ norm for   time inhomogeneous    semigroups.

\section{Heat kernel estimates}

 We first recall a result of \cite{MPZ}. Let
 $$a: [0,T]\times \R^d\to \R^d\otimes\R^d,\ \ b: [0,T]\times\R^d\to \R^d.$$
We consider  heat kernel estimates for the   time dependent second order differential operator
$$L_t^{a,b}:= \ff 1 2 {\rm tr}\{a_t\nn^2\}+\nn_{b_t}$$
 satisfying the following conditions.
 \beg{enumerate} \item[$(H^{a,b})$]  $a_t(x)$ is invertible and there exist constants $C>0$ and $\aa\in (0,1)$ such that
\beg{align*}&\|b_\cdot(0)\|_\infty+\|a\|_\infty+\|a^{-1}\|_\infty\le C,\\
& \sup_{t\in [0,T]} \|a_t(x)-a_t(y)\|\le C|x-y|^\aa,\\
& \sup_{t\in [0,T]} |b_t(x)-b_t(y)|\le C \big(|x-y|+ |x-y|^\aa\big),\ \ \ x,y\in\R^d.\end{align*}
 \item[$(H^{a})$]  $a_t(x)$ is differentiable in $x$,  and there exist constants $C\in (0,\infty) $ and $\aa\in (0,1)$ such that
 $$\|\nn a\|_\infty\le C,\ \ \sup_{t\in [0,T]}\|\nn a_t(x)-\nn a_t(y)\|\le C|x-y|^\aa,\ \ x,y\in \R^d.$$
  \end{enumerate}

Under   $(H^{a,b})$,    for any $s\in [0,T)$, the SDE
$$\d   X_{s,t}^x = b_s( X_{s,t}^x) \d s +\ss{a_s}(  X_{s,t}^x) \d W_s,\ \ t\in [s,T], X_{s,s}^x=x\in\R^d$$
is weakly well-posed with semigroup $\{P_{s,t}^{a,b} \}_{0\le s<t\le T}$ and transition density $\{p_{s,t}^{a,b}\}_{0\le s<t\le T}$ given by
$$P_{s,t}^{a,b} f(x) = \int_{\R^d} p_{s,t}^{a,b}(x,y) f(y)\d y= \E[f(X_{s,t}^x)],\ \ f\in \B_b(\R^d),$$
and  we have the following Kolmogorov backward equation (see Remark 2.2 in \cite{MPZ})
\beq\label{MP2} \pp_s P_{s,t}^{a,b} f=-L_sP_{s,t}^{a,b}f,\ \ f\in C_b^\infty(\R^d), s\in [0,t], t\in (0,T].\end{equation}
 Next, we denote $\psi_{s,t}= \theta_{t,s}^{(1)}$ presented in \cite{MPZ}. Then $(\psi_{s,t})_{0\le s\le t\le T}$ is a family of  diffeomorphisms on $\R^d$ satisfying
\beq\label{MP1} \sup_{0\le s\le t\le T} \big\{\|\nn\psi_{s,t}\|_\infty+   \|\nn \psi_{s,t}^{-1}\|_\infty \big\}\le \dd \end{equation}
for some constant $\dd>0$ depending on $\aa, C$.
For any $\kk>0$, consider the Gaussian heat kernel
 $$p_t^\kk(x):= (\kk \pi t)^{-\ff d 2} \e^{-\ff{|x|^2}{\kk t}},\ \ t>0, \ x\in \R^d.$$
The following result is taken from     \cite[Theorem 1.2]{MPZ}.

\beg{thm}[\cite{MPZ}]\label{M1} Assume $(H^{a,b})$. Then there exist constants $c,\kk>0$ depending on $C,\aa$ such that
\beq\label{MP3}\beg{split} & |\nn^i p_{s,t}^{a,b}(\cdot,y)(x)|\le c (t-s)^{-\ff i 2 } p_{t-s}^\kk(\psi_{s,t}(x)-y),\\
& \ \ \qquad \qquad\qquad i=0,1,2, \ 0\le s<t\le T,\  x,y\in \R^d. \end{split}\end{equation}
If moreover $(H^a)$ holds, then
\beq\label{MP3'} |\nn p_{s,t}^{a,b}(x,\cdot)(y) |\le c (t-s)^{-\ff 1 2} p_{t-s}^\kk(\psi_{s,t}(x)-y),\ \ 0\le s<t\le T, x,y\in\R^d,\end{equation}
and for any $\bb\in (0,1)$ there exists a constant $c' >0$ depending on $C,\aa,\bb$ such that
\beq\label{MP3''}\beg{split} & |\nn p_{s,t}^{a,b}(\cdot,y)(x)-\nn p_{s,t}^{a,b}(\cdot,y')(x) |+ |\nn p_{s,t}^{a,b}(x,\cdot)(y) -\nn p_{s,t}^{a,b}(x,\cdot)(y') |\\
&\le c' |y-y'|^\bb (t-s)^{-\ff {1+\bb} 2} \big\{p_{t-s}^\kk(\psi_{s,t}(x)-y)+ p_{t-s}^\kk(\psi_{s,t}(x)-y')\big\},\\
&\qquad  \  \ 0\le s<t\le T, \ x,x',y\in\R^d.\end{split}\end{equation}
\end{thm}

For any $f\in \B_b(\R^d)\cup \B^+(\R^d)$, let
\beq\label{MMM0} \beg{split}&P_t^\kk f(x):=\int_{\R^d} p_t^\kk(x-y) f(y)\d y,\\
&  \hat P_{s,t}^\kk f(x):=  \int_{\R^d} p_{t-s}^\kk( \psi_{s,t}(x) -y) f(y)\d y,\\
&\tt P_{s,t}^\kk f(x):=   \int_{\R^d} p_{t-s}^\kk( \psi_{s,t}(y) -x) f(y)\d y, \ \ 0\le s<t\le T,\  x\in\R^d.\end{split}\end{equation}
It is well known that for some  constant $c >0$,
\beq\label{M01}\|P_t^\kk\|_{L^p\to L^{p'}}:=\sup_{\|f\|_p\le 1} \|P_t^\kk f\|_{L^{p'}}\le c   t^{-\ff {d(p'-p)}{2pp'}},\ \ t>0, 1\le p\le p'\le \infty.\end{equation}
Combining this with  \eqref{MP1} we obtain
\beq\label{M0} \|\hat P_{s,t}^\kk\|_{L^p\to L^{p'}}+\|\tt P_{s,t}^\kk\|_{L^p\to L^{p'}}\le c   (t-s)^{-\ff {d(p'-p)}{2pp'}},\ \ 0\le s<t\le T, 1\le p\le p'\le \infty.  \end{equation}
for some different constant $c >0$. Below we extend this estimate to the $\tt L^p_q$-$\tt L_q^{p'}$ norm. For any $t\in (0,T]$, let
$$\|f\|_{\tt L_q^p(t)}:= \sup_{z\in \R^d} \bigg(\int_0^t\|1_{B(z,1)} f_s\|_{\tt L^p}^q\d s\bigg)^{\ff 1 q},\ \ p,q\in [1,\infty].$$

\beg{lem}\label{LA} There exists a constant $c >0$ such that for any $0\le s<t\le T, \ 1\le p\le p'\le \infty$ and $ q\in [1,\infty],$
\beq\label{LA1}\beg{split}&     \|  \hat P_{\cdot,t}^\kk   f  \|_{\tt L_q^{p'} (t)}   +  \|  \tt P_{\cdot,t}^\kk   f  \|_{\tt L_q^{p'} (t)}
 \le c \big\|  (t-\cdot)^{-\ff{d(p'-p)}{2pp'}}f\big\|_{\tt L_q^p(t)}, \ \ \ f\in \B^+(\R^d),
\end{split} \end{equation}  where and in the sequel, $(t-\cdot)(s):=t-s$ is a function on $[0,t]$, and
\beq\label{LA2} \sup_{z\in\R^d}  \|g\hat P_{s,t}^\kk(1_{B(z,1)}f)\|_{L^1}\le c(t-s)^{-\ff{d(p'-p)}{2pp'}} \|g\|_{\tt L^{\ff{p'}{p'-1}} } \|f\|_{\tt L^p},\ \ f,g\in \B^+(\R^d). \end{equation}
\end{lem}

  \beg{proof} Let $\BB_n:=\{v\in \mathbb Z^d: |v|_1:=\sum_{i=1}^d |v_i|= n\}, n\ge 0.$
   By  \eqref{MP1},   we find  a constant  $\vv\in (0,1)$ such that for any $n\ge 0$ and $0\le s<t\le T,$
  $$ |\psi_{s,t}(x)-y|^2\ge \vv n^2,\ \ x\in B(\psi_{s,t}^{-1}(z),\vv),\ \ y
  \in \cup_{v\in \BB_n} B(z+v,d),\  z\in\R^d.$$
 Combining  this with \eqref{M01}, we find   constants $c_2,c_3, c_4  >0$ such that for any $z\in \R^d,$   $0\le s<t\le T,$ and $f,g\in \B_b^+(\R^d)$,
  \beg{align*} &\| 1_{B(\psi_{s,t}^{-1}(z),\vv)}g\hat P_{s,t}^\kk   f \|_{L^1}\le \sum_{n=0}^\infty \sum_{v\in \Z^d: |v|_1=n} \|1_{B(\psi_{s,t}^{-1}(z),1)} g\hat P_{s,t}^\kk (1_{B(z+v,d)} f)\|_{L^1}\\
  &\le \sum_{n=0}^\infty \sum_{v\in \BB_n} \int_{\R^d \times \R^d} |1_{B(\psi_{s,t}^{-1}(z),\vv)}g|(x) p_{t-s}^\kk(\psi_{s,t}(x)-y) |1_{B(z+v,d)} f|(y) \d x\d y\\
 &\le c_2 \sum_{n=0}^\infty  \sum_{v\in \BB_n}   \e^{-  \ff{n^2}{c_3(t-s)} }
    \int_{\R^d \times \R^d} |1_{B(\psi_{s,t}^{-1}(z),\vv)}g|(x) p_{2(t-s)}^\kk(\psi_{s,t}(x)- y)  |1_{B(z+v,d)}f|(y)\d x\d y\\
 &\le c_3  \sum_{n=0}^\infty    \sum_{v\in \BB_n}  \e^{-  \ff{n^2}{c_3(t-s)} } \|\{ P_{2(t-s)}^\kk(1_{B(\psi_{s,t}^{-1}(z),\vv)} g ) \}  1_{B(z+v,d)} f\|_{L^1} \\
 &\le c_3  \sum_{n=0}^\infty   \sum_{v\in \BB_n}   \e^{-  \ff{n^2}{c_3(t-s)} } \| P_{2(t-s)}^\kk(1_{B(\psi_{s,t}^{-1}(z),\vv)} g)\|_{L^{\ff p {p-1}} } \| 1_{B(z+v,d)} f\|_{L^p} \\
 &\le c_4 (t-s)^{  -\ff{d(p'-p)}{2pp'}} \|   g \|_{L^{\ff {p'} {p'-1}}}\sum_{n=0}^\infty    \sum_{v\in \BB_n}   \e^{-  \ff{n^2}{c_3(t-s)} }\| 1_{B(z+v,d)} f\|_{L^p}.\end{align*}
Since
\beq\label{LA0} \sup_{z\in \R^d} \bigg(\int_0^t \|1_{B(z, d)}  f \|_{L^p}^q\d s\bigg)^{\ff 1 q} \le c_5 \|  f\|_{\tt L_q^p(t)}\end{equation}
holds for some constant $c_5>0$,  we find a constant $c_6>0$ such that this and H\"older's inequality  imply
\beg{align*}    &\sup_{z\in \R^d} \bigg(\int_0^t  \|1_{B(z,\vv)}\hat P_{s,t}^\kk f_s\|_{L^{p'}}^q \d s\bigg)^{\ff 1 q} =\sup_{z\in \R^d} \bigg(\int_0^t  \|1_{B(\psi_{s,t}^{-1}(z),\vv)}\hat P_{s,t}^\kk f_s\|_{L^{p'}}^q \d s\bigg)^{\ff 1 q}  \\
&\le  \sup_{z\in\R^d} \bigg(\int_0^t \Big\{c_4\|1_{B(z,  d)}  f_s\|_{L^p} (t-s)^{-\ff{d(p'-p)}{2pp'}} \Big\}^q\d s\bigg)^{\ff 1 q}  \sup_{r\in (0,T]}  \sum_{n=0}^\infty    \sum_{v\in \BB_n}   \e^{-  \ff{n^2}{c_3r} }\\
&\le c_6  \big\|(t-\cdot)^{  -\ff{d(p'-p)}{2pp'}} f\big\|_{\tt L_q^p(t)}     \sum_{n=0}^\infty    \sum_{v\in \BB_n}   \e^{-  \ff{n^2}{c_3T} }.\end{align*}
This implies the upper bound for $\hat P^\kk$ in \eqref{LA1},  by noting that  for some constant $K>0$
\beq\label{LA4}   \sum_{n=0}^\infty \sum_{v\in \BB_n}     \e^{-\ff{  n^2}{c_3T}}\le \sum_{n=0}^\infty K(1+n^{d-1}) \e^{-\ff {n^2}{c_3T}}<\infty.\end{equation}
By \eqref{MP1} and integral transforms,   the estimate on $\tt P_{s,t}$  follows  from that of $\hat P_{s,t}^\kk$.

 Similarly, we find a constant $K >1$ such that
 \beg{align*} &\|g\hat P_{s,t}^\kk(1_{B(\psi_{s,t}(z),1)}f)\|_{L^1}\le  \sum_{n=0}^\infty \sum_{v\in \Z^d: |v|_1=n} \| 1_{B(z+v,d)} g\hat P_{s,t}^\kk (1_{B(\psi_{s,t}(z),\vv)} f)\|_{L^1}\\
 &\le \sum_{n=0}^\infty \sum_{v\in \BB_n} \int_{\R^d \times \R^d} |1_{B(z+v,d)}g|(x) p_{t-s}^\kk(\psi_{s,t}(x)-y) |1_{B(\psi_{s,t}(z),\vv)} f|(y) \d x\d y\\
 & \le K   (t-s)^{ -\ff{d(p'-p)}{2pp'}} \sum_{n=0}^\infty
  \sum_{v\in \BB_n}   \e^{-  \ff{n^2}{K(t-s)} } \|g\|_{\tt L^{\ff{p'}{p'-1}}}\| 1_{B(z,\vv)} f\|_{L^p}.\end{align*}
 This together with \eqref{LA0} and \eqref{LA4} implies \eqref{LA2} for some $c>0$.
 \end{proof}

  \section{Proof  of Theorem  \ref{T1}  }

  We first prove assertion $1)$, i.e. the well-posedness of \eqref{E2}.
   For $\gg\in \tt{\scr P}^k_{\nu,T}$, we denote
  \beg{align*}&\si_t^\gg(x):=\si_t(x,\gg_t),\ \ b_t^{\gg,0}(x):= b_t^{(0)}(x,\gg_t(x),\gg_t),\\
 & b_t^\gg(x):= b_t(x,\gg_t(x),\gg_t)= b_t^{(1)}(x)+ b_t^{\gg,0}(x),\ \ \ t\in [0,T], x\in\R^d.\end{align*}

 \beg{lem}\label{L1} Assume {\bf (A)} with $(A_1)$ holding for $\si^\gg$ replacing $\si$ uniformly in $\gg\in \tt{\scr P}^k_{\nu,T}$, where  $k\in [\ff{p_0}{p_0-1},\infty]$.  Then $\eqref{E2}$ is weakly   well-posed for any  $\L_{X_0}\in \tt L^k$ and   $\gg\in \tt{\scr P}^k_{\nu,T},$
and  for any $\bb\in (0,1)$ there exists a constant $c>1$ independent of $\nu$ and $\gg$ such that $\Phi^{\nu}_t\gg:=\ell_{X_t^\gg}$ for $\L_{X_0}=\nu$ satisfies
\beq\label{ETT} \|\Phi^{\nu}\gg\|_{\tt L^k_\infty}\le c \|\nu\|_{\tt L^k}.\end{equation}
Moreover, under the assumption with $(A_2')$ replacing $(A_2)$, the assertion  holds for $({\scr P}^k_{\nu,T}, L^k)$ replacing $(\tt{\scr P}^k_{\nu,T},\tt L^k)$.
\end{lem}
\beg{proof} (a) By $(A_2)$, we have
  \beq\label{MP4} \sup_{\gg\in \tt{\scr P}^k_{\nu,T}} \|b^{\gg,0}\|\le f_0,\ \  \ \  \|f_0\|_{\tt L^{p_0}_{q_0}} <\infty.\end{equation}
According to \cite{YZ}, see also \cite[Theorem 1.1(1)]{W21e},   this together with  $(A_1)$ and $ (A_3)$ imply the well-posedness of \eqref{E2}.
Moreover, by Theorem 6.2.7(ii)-(iii) in \cite{BKRS}, the distribution density function $\ell_{X_t^\gg}$ exists.

(b) To  estimate $ \Phi_t^{\nu} \gg $ for $\gg\in \tt{\scr P}^k_{\nu,T}$, consider the SDE
\beq\label{001} \d \bar X_s^\gg= b_s^{(1)} (\bar X_s^\gg)\d s + \si_s^\gg(\bar X_s^\gg)\d W_s,\ \ s\in [0,t], \ \bar X_0^\gg=X_0^\gg=X_0\ \text{with}\ \L_{X_0}=\nu.\end{equation}
Let   $a^\gg:=\si^\gg(\si^\gg)^*$. Then
$$\E[f(\bar X_t^\gg)]= \E[(P_{0,t}^{a^\gg, b^{(1)} } f)(X_0)] = \int_{\R^d\times \R^d}   p_{0,t}^{a^\gg, b^{(1)} }(x,y) f(y)\nu(\d x)\d y,\ \ f\in \B^+(\R^d),$$
and    \eqref{MP3} holds for $p_{s,t}^{a^\gg, b^{(1)} }$ with constants $c,\kk>0$ uniformly in $\gg$.
So, we find a constant $c_1>0$ such that
\beq\label{**1} \E[f(\bar X_t^\gg)] \le c_1 \int_{\R^d}   (\hat P_{0,t}^\kk f)(x) \nu(\d x)=c_1 (\hat P_{0,t}^{\kk*} \nu)(f),\ \ f\in \B^+(\R^d),\end{equation}
 where
\beq\label{W*}  (\hat P_{0,t}^{\kk *}\nu)(\d y):= \bigg(\int_{\R^d} \hat p_{0,t}^\kk(x,y)\nu(\d x) \bigg)\d y,\ \ t\in (0,T], \nu\in \scr P.\end{equation}
On the other hand, let
$$R_t:= \e^{\int_0^t \<\xi_s,\d W_s\>-\ff 1 2 \int_0^t|\xi_s|^2\d s},\ \ \xi_s:= \big\{\si_s^\gg(\si_s(\si_s^\gg)^*)^{-1}b_s^{\gg,0}\big\}(\bar X_s).$$
By \eqref{MP4}, the uniform boundedness of   $\|\si^\gg(\si^\gg(\si^\gg)^*)^{-1}\|_\infty$, and   Khasminskii's estimate implied by the Krylov's estimate in \cite[Theorem 3.1]{YZ} (see   the proof of \cite[Lemma 4.1(ii)]{XXZZ}), we find a map $K_\gg: [1,\infty)\to (0,\infty)$ such that
\beq\label{KG}  K_\gg(p):= (\E [R_t^p])^{\ff 1 p}<\infty,\ \ p\ge 1.\end{equation}
By Girsanov's theorem,
$$\tt W_s:=W_s-\int_0^s\xi_r\d r,\ \ s\in [0,t]$$
is an $m$-dimensional Brownian motion under the probability measure $\Q_t:=R_t\P$, with which  the SDE \eqref{001} reduces to
$$\d \bar X_s= b_s^\gg(\bar X_s) \d s +\si_s^\gg(\bar X_s)\d \tt W_s,\ \ s\in [0,t], \bar X_0=X_0^\gg.$$
By the weak uniqueness, the law of $X_t^\gg$ under $\P$ coincides with that of $\bar X_t$ under $\Q_t$. Combining this with   \eqref{**1},  \eqref{KG}   and  \eqref{LA2}, for any
$p>1$ and $k'\ge k$ we find  constants $c_1(p), c_2(p)>0$ such that
\beg{align*} & \int_{\R^d} \big\{(\Phi_t^\nu\gg)1_{B(z,1)} f\big\}(y)\d y = \E\big[(1_{B(z,1)} f)(X_t^\gg)\big]=\E\big[R_t(1_{B(z,1)} f)(\bar X_t^\gg)\big]\\
&\le \big(\E[R_t^{\ff{p}{p-1}}]\big)^{\ff {p-1}p} \big(\E\big[(1_{B(z,1)} f^p)(\bar X_t^\gg)\big]\big)^{\ff 1 p}
 \le c_1(p) \bigg(\int_{\R^d}\big\{   (\hat P_{0,t}^\kk (1_{B(z,1)}f^p)\big\}(x)\nu(\d x)\bigg)^{\ff 1 p}\\
&\le c_2(p)\|\nu\|_{\tt L^k}^{\ff 1 p}t^{-\ff{d(k'-k)}{2kk'p}} \|f\|_{\tt L^{\ff {pk'}{k'-1}}},\ \ t\in (0,T], f\in \B^+(\R^d).\end{align*}
Therefore, for any   $\nu\in \tt {\scr P}^k$,
\beq\label{**2}\|\Phi_t^\nu\gg\|_{\tt L^{\ff{pk'}{pk'-k'+1} }}\le c_2 (p)\|\nu\|_{\tt L^k}^{\ff 1 p}t^{-\ff{d(k'-k)}{2kk'p}},\ \ p>1, k'\ge k, \gg\in \tt{\scr P}^k_{\nu,T}, t\in (0,T],\end{equation}
where for $k'=k=\infty$ we set $\ff{pk'}{pk'-k'+1}:=\ff p{p-1},\ \ff{d(k'-k)}{2kk'p}:=0.$ Using \eqref{M0} replacing the estimate in Lemma \ref{LA}, we find a ma $c: (1,\infty)\to (0,\infty)$
such that
\beq\label{**2'}\|\Phi_t^\nu\gg\|_{L^{\ff{pk'}{pk'-k'+1} }}\le c (p)\|\nu\|_{L^k}^{\ff 1 p}t^{-\ff{d(k'-k)}{2kk'p}},\ \ p>1, k'\ge k, \gg\in  L_\infty^k\cap \D_+^1, t\in (0,T].\end{equation}

(c) By the backward Kolmogorov equation \eqref{MP2} and It\^o's formula, for any $f\in C_0^\infty(\R^d)$ we have
\beg{align*} \d \big\{(P_{s,t}^{a^\gg, b^{(1)} }f)(X_s^\gg) \big\}&= \big\{\big(\pp_s + L_s^{a^\gg, b^{(1)} }+\nn_{b_s^{\gg,0}}\big) P_{s,t}^{a^\gg, b^{(1)} }f \big\} (X_s^\gg) \d s + \d M_s\\
&= \big\{\nn_{b_s^{\gg,0}}  P_{s,t}f \big\} (X_s^\gg) \d s + \d M_s,\ \ s\in [0,t]\end{align*}
for some martingale $M_s$. Then
\beq\label{*1}  \beg{split} &\E[f(X_t^\gg)]= \E[P_{t,t}^{a^\gg, b^{(1)} } f(X_t^\gg)]\\
& = \E[P_{0,t}^{a^\gg, b^{(1)} } f(X_0)]+\int_0^t \E\big[(\nn_{b_s^{\gg,0}} P_{s,t}^{a^\gg, b^{(1)} } f)(X_s^\gg)\big]\d s,\ \ s\in [0,t].\end{split} \end{equation}
We explain  that the last term in \eqref{*1} exists. Indeed, by  \cite[Theorem 1.1(2)]{W21e}, there exists a constant $c_2>0$ such that
$$\|\nn P_{s,t}^{a^\gg, b^{(1)} }f\|_\infty \le c_2\|\nn f\|_\infty,\ \ 0\le s\le t, f\in C_b^1(\R^d),$$ so that   \eqref{MP4} and Krylov's estimate (see Theorem 3.1 in \cite{YZ}) yield
  $$\E\bigg(\int_0^t \big|\big(\nn_{b_s^{\gg,0}} P_{s,t}^{a^\gg, b^{(1)} } f\big)(X_s^\gg)\big|\d s\bigg)^{n}\le \E\bigg(\int_0^t c_2\|\nn f\|_\infty |b_s^{\gg,0}|(X_s^\gg)\d s\bigg)^n<\infty,\ \ n\ge 1.$$
Noting that $\Phi_s^{\nu} \gg:= \ell_{X_s^\gg}$ and $P_{s,t}^{a^\gg, b^{(1)} }f(x)= \int_{\R^d} p_{s,t}^{a^\gg, b^{(1)} }(x,y)f(y)\d y$, \eqref{*1}  is equivalent to
\beg{align*}&\int_{\R^d} \big\{\Phi_s^{\nu}f\big\}(y)\d y = \int_{\R^d\times\R^d}  \nu (x) p_{0,t}^{a^\gg, b^{(1)} }(x,y) f(y) \d x\d y \\
&+ \int_0^t \d s \int_{\R^d\times\R^d} (\Phi_s^{\nu}\gg)(x)
\big\{\nn_{b_s^{\gg,0}} p_{s,t}^{a^\gg, b^{(1)} } (\cdot,y)(x)\big\}f(y)\d y,\ \ f\in C_0^\infty(\R^d), s\in [0,t].\end{align*}
  Thus,
\beq\label{*2}\beg{split}  (\Phi_t^{\nu}\gg)(y)= &\,\int_{\R^d}   p_{0,t}^{a^\gg, b^{(1)} }(x,y)\nu(\d x) \\
&+\int_0^t \d s \int_{\R^d} (\Phi_s^{\nu}\gg)(x)  \big\{\nn_{b_s^{\gg,0}} p_{s,t}^{a^\gg, b^{(1)} } (\cdot,y)(x)\big\}\d x,\ \
t\in [0,T].\end{split}\end{equation}
By \eqref{M0} for $p=p'$, $ \|\hat P_t^{\kk*}\nu\|_{\tt L^l} \le K\|\nu\|_{\tt L^l}$ holds for some constant $K>0$.
Combining this with \eqref{MP3},  \eqref{MP4} and \eqref{*2}, we find a constant $c_3>0$ such that
\beq\label{MP5} \|  \Phi_t^{\nu}\gg\|_{\tt L^l} \le c_3 \| \nu\|_{\tt L^l} + c_3\sup_{z\in\R^d} \int_0^t (t-s)^{-\ff 1 2} \big\|1_{B(z,1)} \hat P_{t-s}^\kk \big\{(\Phi_s^{\nu}\gg)
 f_0(s,\cdot)  \big\}\big\|_{L^l} \d s,\ \ l\in [1,\infty].\end{equation}
By $k>k_0  $ and $k\ge \ff{p_0}{p_0-1},$ for any $l\in (k_0,k]\cap [\ff{p_0}{p_0-1},k]$  we have
\beq\label{MP}   q_l:=  \ff{p_0l}{p_0+l}\in (1,l], \ \ \ \ \ff 1 {q_l}=\ff 1 {p_0}+\ff 1 {l},\end{equation}
and   $(p_0,q_0)\in \scr K$ implies
\beq\label{MMP}    \ff 1 2 +\ff{d(l-q_l)}{2lq_l}= \ff 1 2 +\ff d {2 p_0}=:\dd'<\ff {q_0-1}{q_0}.     \end{equation}
 Combining these with  \eqref{LA1}  for $(p',p)=(l, q_l)$ and applying H\"older's inequality,
we find  a  constant   $c_4 >0$ such that
\beg{align*} &\int_0^t (t-s)^{-\ff 1 2} \big\|   \hat P_{t-s}^\kk \big\{(\Phi_s^{\nu}\gg) f_0(s,\cdot)  \big\}\big\|_{\tt L^l}\d s
 \le  c_4  \big\|(t-\cdot)^{-\dd'} f_0 \Phi_\cdot^\nu \gg\big\|_{\tt L_1^{q_l}(t)} \\
&\le c_4  \|f_0\|_{\tt L_{q_0}^{p_0}} \big\|(t-\cdot)^{-\dd'} \Phi_\cdot^\nu\gg\big\|_{\tt L_{\ff{q_0}{q_0-1}}^l(t)},   \ \ l\in   (k_0,k]\cap \Big[\ff{p_0}{p_0-1},k\Big],\end{align*}
where $\{(t-\cdot)^{-\dd'} \Phi_\cdot^\nu\}(s,x):=  (t-s)^{-\dd'} \Phi_s^\nu(x).$ This together with \eqref{MP5} implies that for some constant $c_5>0$,
\beq\label{UP}  \beg{split}   &\|  \Phi_t^{\nu}\gg\|_{\tt L^l}
 \le    c_5  \| \nu\|_{\tt L^l} + c_5  \|f_0\|_{\tt L_{q_0}^{p_0}} \bigg(\int_0^t \Big\{(t-s)^{-\dd'}   \|\Phi_s^{\nu}\gg\|_{\tt L^l}\Big\}^{\ff{q_0}{q_0-1}}  \d s\bigg)^{\ff {q_0-1}{q_0}},\\
&\qquad  \ t\in [0,T], l\in (k_0,k]\cap \Big[\ff{p_0}{p_0-1},k\Big].\end{split} \end{equation}
Similarly, using \eqref{M0} replacing Lemma \ref{LA}, we derive the same estimate for $L^l$ replacing $\tt L^l$ and $\|f_0\|_{  L_{q_0}^{p_0}}$ replacing $\|f_0\|_{\hat L_{q_0}^{p_0}}:$
\beq\label{UP'}  \beg{split}   & \|  \Phi_t^{\nu}\gg\|_{  L^l}  \le     c_5  \| \nu\|_{ L^l} + c_5  \|f_0\|_{  L_{q_0}^{p_0}} \bigg(\int_0^t   \Big\{(t-s)^{-\dd'}   \|\Phi_s^{\nu}\gg\|_{\tt L^l}\Big\}^{\ff{q_0}{q_0-1}} \d s\bigg)^{\ff {q_0-1}{q_0}},\\
&\qquad  \ t\in [0,T], l\in (k_0,k]\cap \Big[\ff{p_0}{p_0-1},k\Big].\end{split} \end{equation}
Below we prove \eqref{ETT} by considering two different situations.

$(c_1)$ $k<\infty$. For any $k'\in (k,\infty)$ we have
$$p_{k,k'}:= \ff{k(k'-1)}{k'(k-1)}>1,\ \ \ff{p_{k,k'}k'}{p_{k,k'}k'-k'+1} =k.$$
Noting that
$$\lim_{k'\downarrow k} \ff{d(k'-k)}{2kk'p_{k,k'}}=0,$$
by \eqref{MMP} we find $k'>k$ such that
$$\vv_{k,k'}:= \ff{d(k'-k)}{2kk'p_{k,k'}}\in \Big(0, 1 -\ff{\dd'q_0}{q_0-1}\Big).$$
Combining this with \eqref{**2} and \eqref{UP} for $l=k$,   we find  a  constant $K>0$ such that
$$\sup_{t\in [0,T]}  \|  \Phi_t^{\nu}\gg\|_{\tt L^k}  \le   K \|\nu\|_{\tt L^k} + K  \sup_{t\in [0,T]}\bigg(\int_0^t (t-s)^{-\ff{q_0\dd'}{q_0-1}}  s^{-\vv_{k,k'}}  \|\nu\|_{\tt L^k}^{\ff{q_0}{p_{k,k'}(q_0-1)}} \d s\bigg)^{\ff {q_0-1}{q_0}}<\infty.$$
Therefore,  by    the generalized
Gronwall inequality (see \cite{YYY}),   \eqref{MMP} and \eqref{UP}   implies  \eqref{ETT}.

When $f_0\in L_{q_0}^{p_0}$, by using   \eqref{**2'} and \eqref{UP'}  replacing \eqref{**2} and \eqref{UP}, we obtain  this estimate  for $L$ replacing $\tt L$.

$(c_2)$ $k=\infty$. We take $k'=k=\infty$, so that
by \eqref{**2}, for any $p>1$ we find a constant $c(p)>0$ such that
  $$\|\Phi_t^\nu\gg\|_{\tt L^{\ff p{p-1}}}\le c(p) \|\nu\|_{\tt L^k}^{\ff 1 p}.$$
  Combining this with \eqref{UP} for   $l \in (\ff{p_0}{p_0-1}\lor k_0,\infty)$ and $p:=\ff l{l-1}>1$,    we obtain
  $$\sup_{t\in [0,T]}    \|  \Phi_t^{\nu}\gg\|_{\tt L^{l} } <\infty,$$
  so that by the generalized Gronwall inequality, \eqref{UP} implies      \eqref{ETT}  for  $l  \in   (\ff{p_0}{p_0-1}\lor k_0,\infty)$ replacing $k=\infty$ with  a uniform constant $c>0$.
  By letting $l\uparrow k=\infty$, we prove \eqref{ETT}.

 Noting that a probability  density function $\rr\in L^\infty$ implies $\rr\in L^l$ for any $l\ge 1$, when $f_0\in L_{q_0}^{p_0}$
 we prove   \eqref{ETT} for $L$ replacing $\tt L$ by using   \eqref{**2'} and \eqref{UP'}  replacing \eqref{**2} and \eqref{UP}.

  \end{proof}

 \beg{proof}[Proof of Theorem \ref{T1}(1)]    By Lemma \ref{LA}, \eqref{E2} is weakly well-posed. By \cite{YZ} or \cite{W21e}, it is also strongly well-posed provided \eqref{SI} holds.
 Thus, as explained in the end of Section 1 that for the weak or strong well-posedness of \eqref{E1}, it suffices to prove that $\Phi^\nu$ has a unique fixed point in $\tt{\scr P}^k_{\nu,T}$. In general, for   any $\nu_1,\nu_2\in \tt {\scr P}^{k}$ and $\gg^1,\gg^2\in \tt{\scr P}^k_{\nu,T}$, we estimate
 $$\tt d_{k,\ll}(\Phi^{\nu_1
 }\gg^1,\Phi^{\nu_2}\gg^2):=\sup_{t\in [0,T]} \e^{-\ll t} \|\Phi^{\nu_1}_t\gg^1
 -\Phi^{\nu_2}_t\gg^2\|_{\tt L^k},\ \  \ll>0.$$

By \eqref{*2},  $(A_2)$ and  \eqref{MP3}, we find a constant $c_1>0$ such that
\beq\label{GG1} \beg{split}& \|  \Phi_t^{\nu_1}\gg^1- \Phi_t^{\nu_2}\gg^2 \|_{\tt L^k}- c_1\|\nu_1-\nu_2\|_{\tt L^k} \\
 &\le c_1 \int_0^t (t-s)^{-\ff 1 2}  \Big\|  \hat P_{s,t}^\kk\Big\{f_0(s,\cdot)\big[|\Phi_s^{\nu_1}\gg^1- \Phi_s^{\nu_2}\gg^2|\\
 &\qquad\qquad \qquad + s^\theta(\Phi_s^{\nu_1}\gg^1)\big( |\gg_s^1-\gg_s^2| +\|\gg_s^1-\gg_s^2\|_{\tt L^k}\big)\big]\Big\}\Big\|_{\tt L^{k}}\d s.\end{split}\end{equation}
 Letting
 $$F_l(s,x):= (t-s)^{-\ff{d(k-l)}{2kl} -\ff 1 2}  \big[|\Phi_s^{\nu_1}\gg^1- \Phi_s^{\nu_2}\gg^2|+ s^\theta(\Phi_s^{\nu_1}\gg^1)\big( |\gg_s^1-\gg_s^2| +\|\gg_s^1-\gg_s^2\|_{\tt L^k}\big)\big](x) $$ for $l\in [1,\ff{kp_0}{k+p_0}],$
 by \eqref{LA1} for $q=1$ and $(p',p)=(k,l)$, and applying H\"older's inequality,  we find a constant $c_2>0$ such that
 \beg{align*} &\int_0^t (t-s)^{-\ff 1 2}  \Big\|  \hat P_{s,t}^\kk\Big\{f_0(s,\cdot)\big[|\Phi_s^{\nu_1}\gg^1- \Phi_s^{\nu_2}\gg^2|+ s^\theta(\Phi_s^{\nu_1}\gg^1)\big( |\gg_s^1-\gg_s^2| +\|\gg_s^1-\gg_s^2\|_{\tt L^k}\big)\big]\Big\}\Big\|_{\tt L^k}\d s\\
 &\le c_2\| f_0 F_l\|_{\tt L_1^l}\le c \|f_0\|_{\tt L_{q_0}^{p_0} } \|F_l\|_{\tt L_{\ff{q_0}{q_0-1}}^{\ff{p_0l}{p_0-l}}}\\
 &\le c_2 \|f_0\|_{\tt L_{q_0}^{p_0} }  \bigg(\int_0^t \Big\{  (t-\cdot)^{-\ff{d(k-l)}{2kl} -\ff 1 2}   \Big[\big\|\Phi^{\nu_1}\gg^1- \Phi^{\nu_2}\gg^2\big\|_{\tt L^{\ff{p_0l}{p_0-l}}}\\
 &\qquad\qquad\qquad \qquad + s^\theta\big\|\big(\Phi_s^{\nu_1}\gg^1\big)\big( |\gg_s^1-\gg_s^2| +\|\gg_s^1-\gg_s^2\|_{\tt L^k}\big)\big\|_{\tt L^{\ff{p_0k}{p_0-k}}}  \Big] \Big\}^{\ff {q_0}{q_0-1}} \d s\bigg)^{\ff{q_0-1}{q_0}}. \end{align*}  Since $l\in [1,\ff{kp_0}{k+p_0}] $ implies
 $ \ff{p_0l}{p_0-l}\le k,$ combining this with \eqref{GG1} and applying H\"older's inequality,  we find a constant $c_3>0$ such that
\beq\label{GG5}\beg{split} & \|  \Phi_t^{\nu_1}\gg^1- \Phi_t^{\nu_2}\gg^2 \|_{\tt L^k}- c_1\|\nu_1-\nu_2\|_{\tt L^k} \\
&\le c_3  \bigg(\int_0^t \Big\{ (t-s)^{-\ff{d(k-l)}{2kl} -\ff 1 2}   \Big[\big\|\Phi_s^{\nu_1}\gg^1- \Phi_s^{\nu_2}\gg^2\big\|_{\tt L^{k}}\\
 &\qquad  + s^\theta \big\| \Phi_s^{\nu_1}\gg^1\big\|_{\tt L^{\ff{kp_0l}{k(p_0-l)-p_0l}} } \big\|\gg_s^1-\gg_s^2\big\|_{\tt L^k} \Big] \Big\}^{\ff {q_0}{q_0-1}} \d s\bigg)^{\ff{q_0-1}{q_0}},\ \ l\in \Big[1,  \ff{kp_0}{k+p_0}\Big].\end{split} \end{equation}  Letting
\beq\label{GG2} \aa_l:= \ff{q_0 }{q_0-1} \Big( \ff{d(k-l)}{2kl} + \ff 1 2\Big),\ \ \bb_l:=\ff{kp_0l}{k(p_0-l)-p_0l},\end{equation}
  by the definition of $\tt d_{k,\ll}$, this implies that for any $\ll>0$ and $l\in  [1,\ff{kp_0}{k+p_0}],$
 \beq\label{WM} \beg{split} & \tt d_{k,\ll}(\Phi^{\nu_1}\gg^1, \Phi^{\nu_2}\gg^2) \le   c_1\|\nu_1-\nu_2\|_{\tt L^k} \\
 &+ c_3\tt d_{k,\ll}(\Phi^{\nu_1}\gg^1, \Phi^{\nu_2}\gg^2)
   \sup_{t\in (0,T]}\bigg(\int_0^t (t-s)^{- \aa_l}  \e^{-\ff{\ll q_0}{q_0-1}(t-s)} \d s\bigg)^{\ff {q_0-1}{q_0}}\\
 & + c_3\tt d_{k,\ll}(\gg^1,\gg^2) \sup_{t\in (0,T]}\bigg(\int_0^t  (t-s)^{- \aa_l}  \e^{-\ff{\ll q_0}{q_0-1}(t-s)}
\big(s^\theta  \|\Phi_s^{\nu_1}\gg^1 \|_{\tt L^{\bb_l}}\big)^{\ff{q_0}{q_0-1}}\d s\bigg)^{\ff {q_0-1}{q_0}}.\end{split}\end{equation}

 Below we complete the proof  by considering two different situations respectively.

 \

 {\bf  (a)}  Let $k<\infty$.   By $(p_0, q_0)\in\scr K$ and $k>k_0:=\ff{d}{2\theta+1-dp_0^{-1}-2q_0^{-1}}$, $\aa_l$ in \eqref{GG2} satisfies
 $$\lim_{l\uparrow \ff{kp_0}{k+p_0}} \aa_l+ \ff{q_0}{q_0-1}\Big(\ff d {2k} -\theta\Big)^+= \ff{q_0}{q_0-1}\Big\{\ff d{2p_0} +\ff 1 2 +\Big(\ff d {2k}-\theta\Big)^+\Big\}<1.$$
 So, we may take $l\in (1,  \ff{kp_0}{k+p_0})$ such that
\beq\label{GG3} \aa_l+ \ff{q_0}{q_0-1}\Big(\ff d {2k} -\theta \Big)^+<1,\ \ \ \bb_l\in (1,\infty).\end{equation}
   By \eqref{**2} for $k'=\infty$ and $p=\ff{\bb_l}{\bb_l-1}$,  there exists a constant $c_4>0$ such that
 $$\|\Phi_s^{\nu_1}\gg^1 \|_{ L^{\bb_l}}\le c_4 \|\nu_1\|_{\tt L^k} s^{-\ff d{2k}}.$$
 Combining this with \eqref{WM}  and  \eqref{GG3},  when $\ll$ is large enough  increasing in  $\|\nu_1\|_{\tt L^k} (\le\|\nu_2\|_{\tt L^k}$), we obtain
 $$\tt d_{k,\ll}(\Phi^{\nu_1}\gg^1, \Phi^{\nu_2}\gg^2) \le   c_1\|\nu_1-\nu_2\|_{\tt L^k} \\
 + \ff 1 4 \tt d_{k,\ll}(\Phi^{\nu_1}\gg^1, \Phi^{\nu_2}\gg^2)  +\ff 1 4 \tt d_{k,\ll}(\gg^1,\gg^2).$$
 Taking $\nu_1=\nu_2=\nu$ we prove the contraction of $\Phi^\nu$ on the complete metric space  $(\tt{\scr P}^k_{\nu,T}, \tt d_{k,\ll})$, and hence $\Phi^\nu$  has a unique fixed point.
This implies the weak (also strong under \eqref{SI}) well-posedness of \eqref{E1}. Moreover, for two solutions $(X^i)_{i=1,2}$ of this SDE with initial distribution densities $(\nu_i)_{i=1,2}$,
by taking $\gg^i= \L_{X_\cdot^i}$ we have  $\gg^i=\Phi^{\nu_i}\gg^i$, so that this estimate implies \eqref{S2} for some increasing function $\LL$.

\

{\bf (b)} Let  $k=\infty$. By taking $l=p_0$, we have $\bb_l=\infty$ and $\theta>\ff 2 {q_0}+\ff d{p_0}-1$ in $(A_2)$ implies
$$ \aa_l+ \ff{q_0}{q_0-1}\Big(\ff d {2k} -\theta \Big)^+ =\ff{q_0}{q_0-1}\Big\{\ff d{2p_0}+\ff 1 2+ \theta^-  \Big\}<  1.$$
Combining   \eqref{WM} with   \eqref{ETT} for $k=\infty$,  we derive that for a large enough $\ll>0$ increasing in $\|\nu_1\|_{\tt L^\infty}\ (\le  \|\nu_2\|_{\tt L^\infty}$),
\beg{align*} &\tt d_{k,\ll}(\Phi^{\nu_1}\gg^1, \Phi^{\nu_2}\gg^2) \le   c_1\|\nu_1-\nu_2\|_{\tt L^\infty} \\
 &+ c_3\tt d_{k,\ll}(\Phi^{\nu_1}\gg^1, \Phi^{\nu_2}\gg^2)
   \sup_{t\in (0,T]}\bigg(\int_0^t (t-s)^{-\aa_l} \e^{-\ff{\ll q_0}{q_0-1}(t-s)} \d s\bigg)^{\ff {q_0-1}{q_0}}\\
 & + c_3\tt d_{k,\ll}(\gg^1,\gg^2) \sup_{t\in (0,T]}\bigg(\int_0^t  (t-s)^{-\aa_l} \e^{-\ff{\ll q_0}{q_0-1}(t-s)}
(s^\theta  \|\nu_1\|_{\tt L^\infty})^{\ff {q_0}{q_0-1}}\d s\bigg)^{\ff {q_0-1}{q_0}}\\
&\le   c_1\|\nu_1-\nu_2\|_{\tt L^\infty}
  + \ff 1 4 \tt d_{k,\ll}(\Phi^{\nu_1}\gg^1, \Phi^{\nu_2}\gg^2)
    + \ff 1 4 \tt d_{k,\ll}(\gg^1,\gg^2).\end{align*}
Then we finish the proof as shown in step (a).
  \end{proof}

\beg{proof}[Proof of Theorem \ref{T1}(2)]  Let {\bf (A)} hold for $(A_2')$ replacing $(A_2)$. By \eqref{M0} and H\"older's inequality,  we find    constants $c_1,c_2>0$ such that
for any $ 0\le s<t\le T$ and $l\in [1,\ff{kp_0}{k+p_0}],$
\beg{align*} &\Big\|\hat P_{s,t}^\kk\Big\{\big(C+f_0(s,\cdot)\big) \big(|\Phi_s^{\nu_1}\gg^1- \Phi_s^{\nu_2}\gg^2|+ s^\theta(\Phi_s^{\nu_1}\gg^1) |\gg_s^1-\gg_s^2|\big)\Big\} \Big\|_{L^k}\\
&\le  c_1(t-s)^{-\ff {d(k-l)}{2kl}} \Big\{  \big\| \Phi_s^{\nu_1}\gg^1- \Phi_s^{\nu_2}\gg^2\|_{L^l}    +  s^\theta \big\| (\Phi_s^{\nu_1}\gg^1) |\gg_s^1-\gg_s^2|   \big\|_{L^l}\\
&\qquad +\big\|f_0(s,\cdot) (\Phi_s^{\nu_1}\gg^1- \Phi_s^{\nu_2}\gg^2)\big\|_{L^l}+ s^\theta \big\|f_0(s,\cdot)(\Phi_s^{\nu_1}\gg^1) |\gg_s^1-\gg_s^2|   \big\|_{L^l}\Big\}\\
&\le  c_1(t-s)^{-\ff {d(k-l)}{2kl}} \Big\{  \big\| \Phi_s^{\nu_1}\gg^1- \Phi_s^{\nu_2}\gg^2\|_{L^l}   +
s^\theta \big\| \Phi_s^{\nu_1}\gg^1  \big\|_{L^{\ff{kl}{k-l}} }\big\|\gg_s^1-\gg_s^2|   \big\|_{L^k} \\
&\quad +\big\|f_0(s,\cdot)\big\|_{L^{p_0}}  \Big(\big\|\Phi_s^{\nu_1}\gg^1- \Phi_s^{\nu_2}\gg^2 \big\|_{L^{\ff {p_0l}{p_0-l}}}+ s^\theta
\big\|\Phi_s^{\nu_1}\gg^1\big\|_{L^{\ff{p_0kl}{p_0k-kl-p_0l}} } \big\|\gg_s^1-\gg_s^2   \big\|_{L^k}\Big)\Big\}.\end{align*}
 Noting that $l\in [1,\ff{kp_0}{k+p_0}]$ implies $l\lor \ff {p_0l}{p_0-l}\le k$ and $\ff{kl}{k-l}\le \ff{p_0kl}{p_0k-kl-p_0l}$,  by combining this with $(A_2')$, \eqref{MP3}, \eqref{*2} and H\"older's inequality, we find   constants $c_3, c_4>0$ such that
    \beg{align*} & \|  \Phi_t^{\nu_1}\gg^1- \Phi_t^{\nu_2}\gg^2 \|_{L^k}- c_1\|\nu_1-\nu_2\|_{L^k} \\
    &\le c_3 \int_0^t (t-s)^{-\ff 1 2} \Big\|\hat P_{s,t}^\kk\Big\{\big(C+f_0(s,\cdot)\big) \big(|\Phi_s^{\nu_1}\gg^1- \Phi_s^{\nu_2}\gg^2|+ s^\theta(\Phi_s^{\nu_1}\gg^1) |\gg_s^1-\gg_s^2|\big)\Big\} \Big\|_{L^k}\d s\\
 &\le c_4 \big(1+\|f_0\|_{\tt L_{q_0}^{p_0}}\big)  \bigg(\int_0^t \Big\{\|(t-s)^{-\ff{d(k-l)}{2kl} -\ff 1 2}   \Big[\|\Phi_s^{\nu_1}\gg^1- \Phi_s^{\nu_2}\gg^2\|_{L^{k}}\\
 &\qquad\qquad\qquad + s^\theta \| \Phi_s^{\nu_1}\gg^1\|_{L^{\ff{kp_0l}{k(p_0-l)-p_0l}} } |\gg_s^1-\gg_s^2\|_{L^k} \Big] \Big\}^{\ff {q_0}{q_0-1}} \d s\bigg)^{\ff{q_0-1}{q_0}},\ \ l\in \Big[1,  \ff{kp_0}{k+p_0}\Big].\end{align*}
  Then the remainder of the proof is similar to that of Theorem \ref{T1}(1) from \eqref{GG5} with $L$ replacing   $\tt L$.
    \end{proof}

\section{Proof of Theorem \ref{T2}}

Let $\nu\in \scr P^\infty$ with $\ell_{\nu}\in C_b^\bb$. By Theorem \ref{T1} and {\bf (B)},
for any
$\gg\in  \scr P^\infty_{\nu,T}$,    the  following density dependent   SDE has a unique (weak and strong) solution with $\ell_{X_\cdot^{\gg,\nu}}\in L_\infty^\infty$:
\beq\label{GM} \d X_t^{\gg,\nu} = b_t(X_t^{\gg,\nu}, \ell_{X_t^{\gg,\nu}}(X_t^{\gg,\nu}),\ell_{X_t^{\gg,\nu}})\d t + \si_t^{\gg}(X_t^{\gg,\nu}),\ \ \L_{X_0^{\gg,\nu}}=\nu, t\in [0,T],\end{equation}   and there exists a constant $c>0$ depending on $C,\aa$ such that
\beq\label{MC0} \|\ell_{X_t^{\gg,\nu}}\|_\infty\le c \|\ell_\nu\|_\infty,\ \  \gg\in L_\infty^\infty\cap \D_+^1. \end{equation}
We aim to show that the map
$$\gg\mapsto  \ell_{X_\cdot^{\gg,\nu}}$$ has a unique fixed point in $L_\infty^\infty\cap \D_+^1$, such that the (weak and strong) well-posedness of \eqref{GM} implies that of \eqref{E1}.
As shown   in the proof of Theorem \ref{T1},  we will need heat kernel estimates presented in Section 2 for the operator $L_t^{a^\gg, b^{\gg,\nu}}$, where
$$a_t^\gg:= \ff 1 2 \si^\gg_t(\si_t^\gg)^*,\ \ \ b_t^{\gg,\nu}:= b_t\big(\cdot, \ell_{X_t^{\gg,\nu}}(\cdot), \ell_{X_t^{\gg,\nu}}\big),\ \ t\in [0,T].$$
To this end, we first prove the H\"older continuity of $b_t^{\gg,\nu}$. By {\bf (B)}, this follows from the H\"older continuity of $\ell_{X_t^{\gg,\nu}}.$

\beg{lem}\label{L2} Assume {\bf (B)} and let $\bb\in(0, 1-\ff d {p_0}-\ff 2 {q_0})$.  Then there exists a constant $c>0$ such that for any $\gg\in L_\infty^\infty\cap \D_+^1$ and   $\nu\in \scr P^\infty$ with $\ell_{\nu}\in C_b^\bb$,
\beq\label{LB1} \|\ell_{X_t^{\gg,\nu}}\|_{C_b^\bb}\le c \|\ell_\nu\|_{C_b^\bb},\ \ \ t\in (0,T].\end{equation}

\end{lem}
\beg{proof}   Simply denote $\ell_t= \ell_{X_t^{\gg,\nu}}.$ Let $p_{s,t}^\gg$ be the heat kernel for the operator
$$L_t^{\gg}:= \ff 1 2 {\rm div} \big\{a_t^\gg\nn \big\}= L_t^{a^\gg, \bar b^\gg},$$
where
$$a_t^\gg:= \ff 1 2 \si^\gg_t(\si_t^\gg)^*,\ \ \ (\bar b_t^{\gg})_i:=\ff 1 2  \sum_{j=1}^d \pp_j (a_t^\gg)_{ij}.$$
Then $p_{s,t}^\gg(x,y)=p_{s,t}^\gg(y,x)$, and by {\bf (B)} and Theorem \ref{M1}, there exist  constants $c,\kk>0$ depending on $C,\aa,\bb$ such that for some  diffeomorphisms   $\{\psi_{s,t}\}_{0\le s\le t\le T}$ satisfying \eqref{MP1},
\beq\label{MC1} \beg{split} & |\nn^i p_{s,t}^\gg(\cdot,y)(x)|\le c_1 (t-s)^{-\ff i 2} p_{t-s}^\kk (\psi_{s,t}(x)-y),\ \ i=0,1,2,\\
 & |\nn  p_{s,t}^\gg(x,\cdot)(y)|\le c_1 (t-s)^{-\ff 1 2} p_{t-s}^\kk (\psi_{s,t}(x)-y),\\
 &|\nn  p_{s,t}^\gg(\cdot,y)(x)- \nn  p_{s,t}^\gg(\cdot,y')(x)|\\
 &\le c_1|y-y'|^\bb (t-s)^{-\ff {1+\bb} 2}\big\{ p_{t-s}^\kk (\psi_{s,t}(x)-y)+p_{t-s}^\kk (\psi_{s,t}(x)-y')\big\},\\
 &\qquad  \  0\le s<t\le T, \ x,y,y'\in \R^d.\end{split}\end{equation}
 By the argument leading to \eqref{*2} for $\bar b^\gg$ replacing $b^{(1)}$, we obtain
\beq\label{MC} \ell_t(y)= \int_{\R^d} p_{0,t}^\gg(x,y)\ell_{\nu}(x)\d x + \int_0^t\d s \int_{\R^d} \ell_s(x)\big\{\nn_{b_s(x, \ell_s(x),\ell_s)-\bar b_s^\gg(x) } p_{s,t}^\gg(\cdot,y)\big\}(x)\d x.\end{equation}
By the symmetry of $p^\gg_{0,t}(x,y)$ we have
\beq\label{MC2} \int_{\R^d} p_{0,t}^\gg(x,y)\ell_{\nu}(x)\d x = \int_{\R^d} p_{0,t}^\gg(y,x)\ell_{\nu}(x)\d x=: (P_{0,t}^\gg \ell_{\nu})(y).\end{equation}
Let  $X_t^x$ solve the SDE
$$\d X_t^x= \bar b_t^\gg(X_t^x)\d t +\si_t^\gg(X_t)\d W_t,\ \ t\in [0,T], X_0=x.$$
By \cite[(4.8)]{XXZZ},  we find a constant $c_1>0$  depending on $C,\aa$ in  {\bf (B)} such that
$$\E\Big[\sup_{t\in [0,T]} |X_t^x-X_t^y|\Big]\le c_1 |x-y|,\ \ x,y\in \R^d.$$
Then  \eqref{MC2} implies
\beg{align*} &|(P_{0,t}^\gg \ell_{\nu})(y)-(P_{0,t}^\gg \ell_{\nu})(y')| =|\E[\ell_{\nu}(X_t^y)-\ell_{\nu}(X_t^{y'})]|\\
&\le \|\ell_{\nu}\|_{C_b^\bb} \E[|X_t^y-X_t^{y'}|^\bb] \le \|\ell_{\nu}\|_{C_b^\bb} (c_1 |y-y'|)^\bb.\end{align*}
Since {\bf (B)} implies
$|b|+|\bar b^\gg|\le c f_0$  for some constant $c>0$, by combining this with   \eqref{MC0}, the last inequality in \eqref{MC1},  and \eqref{MC},    we find a  constant  $c_2  >0$ independent of $\gg,\nu$ such that
\beg{align*} & |\ell_t(y)-\ell_t(y')|- c_2 |y-y'|^\bb\\
&\le  c_2 \|\ell_{\nu}\|_\infty |y-y'|^\bb  \int_0^t  (t-s)^{-\ff {1+\bb} 2} \big\{\tt P_{s,t}^\kk f_0(s,\cdot)(y)  +\tt P_{s,t}^\kk f_0(s,\cdot) (y') \big\} \d s,
 \end{align*} where $\tt P_{s,t}$ is in \eqref{MMM0}.
By \eqref{LA1} for $(p,q)=(p_0,q_0)$ and $p'=\infty$, we find a constant $c_3>0$ such that this   implies
\beg{align*} &\ff{|\ell_t(y)-\ell_t(y')|- c_2 |y-y'|^\bb}{c_2 \|\ell_{\nu}\|_\infty |y-y'|^\bb}\\
&\le    \int_0^t  (t-s)^{-(\ff {1+\bb} 2+\ff{d}{2p_0})}
  \Big(\tt P_{s,t}^\kk \big\{(t-s)^{\ff{d}{2p_0}} f_0(s,\cdot) \big\}(y)+\tt P_{s,t}^\kk \big\{(t-s)^{\ff{d}{2p_0}}  f_0(s,\cdot) \big\} (y') \Big) \d s\\
&\le 2c_2 \bigg(\int_0^t (t-s)^{-(\ff {1+\bb} 2+\ff{d}{2p_0})\ff{q_0}{q_0-1}} \d s\bigg)^{\ff{q_0-1}{q_0}} \big\|\tt P_{\cdot,t}^\kk\{(t-\cdot)^{\ff d{2p_0}} f_0\}\big\|_{  L_{q_0}^\infty(t)} \\
&\le c_3 \|f\|_{\tt L_{q_0}^{p_0}}, \ \ y\ne y', t\in (0,T],
 \end{align*} where we have used the fact that $\|\cdot\|_{\tt L_{q_0}^\infty}= \|\cdot\|_{L_{q_0}^\infty}$ and $(\ff {1+\bb} 2+\ff{d}{2p_0})\ff{q_0}{q_0-1}<1$ due to $\bb\in (0,1-\ff 2{q_0}-\ff d {p_0})$.  Combining this with \eqref{MC0}, we finish the proof.
\end{proof}

The next lemma contains two classical estimates on  the operator $1 -\DD$ and the heat semigroup $P_t=\e^{t\DD}$.

 \beg{lem}\label{L3} Let $P_t=\e^{t\DD}$.
 \beg{enumerate} \item[$(1)$] For any $\bb> 0$, there exists a constant $c>0$ such that
 $$\|(1-\DD)^{\ff \bb 2} f\|_\infty\le c \|f\|_{C_b^{\bb}}.$$
 \item[$(2)$] For any $\aa,\bb,k\ge 0$, there exists a constant $c>0$ such that
 $$\|(1-\DD)^{-k}P_t f\|_{C_b^{\aa+\bb}}\le c t^{-(\ff \aa 2-k)^+} \|f\|_{C_b^\bb},\ \ t>0.$$\end{enumerate} \end{lem}

\beg{proof}[Proof of Theorem \ref{T2}]   For   $ \L_{X_0^i}=\nu_i$ with $\ell_{\nu_i} \in C_b^\bb(\R^d)$ and $\gg^i\in L_\infty^\infty\cap \D_+^1,$
simply denote $$\ell_t^i=\ell_{X_t^{\gg^i,\nu_i}},\ \  b^{\ell^i}_t:= b_t(\cdot,\ell_t^i(\cdot), \ell_t^i),\ \ t\in [0,T], \ i=1,2.$$ Without loss of generality, let $\|\ell_{\nu_2}\|_{C_b^\bb}\le \|\ell_{\nu_1}\|_{C_b^\bb}.$

By \eqref{MC} with $(\nu,\gg)=(\nu_1, \gg^1)$, we obtain
$$  \ell_t^1(y)= P_{0,t}^{\gg^1} \ell_{\nu_1}(y)  + \int_0^t\d s \int_{\R^d} \ell_s^1(x)\big\{\nn_{b_s^{\ell^1}(x)-\bar b_s^{\gg^1}(x) } p_{s,t}^{\gg^1}(\cdot,y)\big\}(x)\d x.$$
By  the argument leading to \eqref{*2} for  $(p_{s,t}^{\gg^1} , X_s^{\gg^2,\nu_2})$ replacing $(p_{s,t}^{a^\gg, b^{(1)}}, X_s^\gg)$, we derive
  \beg{align*} \ell_t^2(y)= &\,P_{0,t}^{\gg^1} \ell_{\nu_2}(y)
+ \int_0^t\d s \int_{\R^d} \ell_s^2(x)\big\{\nn_{b_s^{\ell^2}(x)-\bar b_s^{\gg^1}(x) } p_{s,t}^{\gg^1}(\cdot,y)\big\}(x)\d x\\
&+\ff 1 2 \sum_{i,j=1}^d \int_0^t\d s \int_{\R^d} \big\{\ell_s^2 (a_s^{\gg^2}-a_s^{\gg^1})_{ij} \pp_i\pp_j  p_{s,t}^{\gg^1}(\cdot,y)\big\}(x)\d x.\end{align*}
Thus,
\beq\label{MC4} \|\ell_t^1 -\ell_t^2 \|_\infty \le I_1+I_2+\sum_{i,j=1}^d I_{ij}, \end{equation}
where
\beq\label{MMC}  I_1:= \|P_{0,t}^{\gg^1}\ell_{\nu_1}- P_{0,t}^{\gg^1} \ell_{\nu_2}\|_\infty\le  \|\ell_{\nu_1}-\ell_{\nu_2}\|_\infty,\end{equation}
and
\beg{align*}&  I_2:=\int_0^t \d s\int_{\R^d} \Big|\Big\{\big[\ell_s^2  ( b_s^{\ell^2} -\bar b_s^{\gg^1}  )- \ell_s^1 ( b_s^{\ell^1} -\bar b_s^{\gg^1} )\big] \nn p_{s,t}^{\gg^1}(\cdot,y)\Big\}(x)\Big|\d x,\\
& I_{ij}:= \ff 1 2  \sup_{y\in\R^d} \bigg|\int_0^t\d s \int_{\R^d} \big\{\ell_s^2  ( a_s^{\gg^2}-a_s^{\gg^1})_{ij} \pp_i\pp_j  p_{s,t}^{\gg^1}(\cdot,y)\big\}(x)\d x\bigg|.\end{align*}
Below we estimate $I_2$ and $I_{ij}$   respectively.

 Firstly, by  {\bf (B)} and \eqref{MC0}, we find a constant $c_1>0$ such that
\beg{align*} &\big|\ell_s^2   \{ b_s^{\ell^2}(x)-\bar b_s^{\gg^1}(x) \}- \ell_s^1(x) \{ b_s^{\ell^1}(x)-\bar b_s^{\gg^1}(x)\}\big|\\
&\le \|\ell_s^1-\ell_s^2\|_\infty |b_s^{\ell^2}(x)-\bar b_s^{\gg^1}(x) | + \|\ell_s^2\|_\infty |b_s^{\ell^2}(x)- b_s^{\ell^1}(x)|\\
&\le c_1\|\ell_{\nu_2}\|_\infty  \|\ell_s^1-\ell_s^2\|_\infty f_0(s,x),\ \ s\in [0,T], x\in\R^d. \end{align*}
 Combining this with
  \eqref{MC1} for $i=1$,  \eqref{LA1} for $(p,q)=(p_0,q_0)$ and $p'=\infty$, and applying H\"older's inequality, we find    constant $c_2, c_3>0$ such that
 \beq\label{NW*1}\beg{split}  & I_2  \le c_2\|\ell_{\nu_2}\|_\infty\int_0^t (t-s)^{-\ff 1 2}  \|\ell_s^1-\ell_s^2\|_\infty \tt P^\kk_{s,t}  f_0(s,\cdot)(y) \d s\\
 &\le   c_2 \|\ell_{\nu_2}\|_\infty \bigg(\int_0^t \big\{(t-s)^{-  (\ff 1  2+\ff d {2p_0}) }   \|\ell_s^1-\ell_s^2\|_\infty\big\}^{\ff{q_0}{q_0-1}}\d s\bigg)^{\ff {q_0-1} {q_0 }}
 \|(t-\cdot)^{\ff d {2p_0}}f_0\|_{\tt L_{q_0}^\infty}\\
 &\le c_3\|\ell_{\nu_2}\|_\infty\|f_0\|_{\tt L_{q_0}^{p_0}}   \bigg(\int_0^t (t-s)^{-\ff{q_0}{q_0-1} (\ff 1  2+\ff d {2p_0}) }   \|\ell_s^1-\ell_s^2\|_\infty^{\ff{q_0}{q_0-1}}\d s\bigg)^{\ff {q_0-1} {q_0 }} ,\ \ t\in [0,T]. \end{split}\end{equation}

  Next, by integration by parts formula, {\bf (B)}, \eqref{LB1},       \eqref{MC1} for $i=1$ and Lemma \ref{L3},   for any $\dd:=\aa\land \bb$, we  find   constants  $c_4, c_5>0$   such that
 \beg{align*} &\bigg| \int_{\R^d} \big\{\ell_s^2 (a_s^{\gg^2}-a_s^{\gg^1})_{ij} \pp_i\pp_j  p_{s,t}^{\gg^1}(\cdot,y)\big\}(x)\d x\bigg| \\
 & =\bigg| \int_{\R^d}  \Big[ (1-\DD)^{\ff \dd 2}  \big\{ \ell_s^2    (a_s^{\gg^2}-a_s^{\gg^1})_{ij} \big\}(x)\Big]\cdot    \Big[\pp_i\pp_j (1-\DD)^{-\ff \dd 2}  p_{s,t}^{\gg^1}(\cdot,y)  (x)\Big] \d x\bigg|\\
 &\le \big\|(1-\DD)^{\ff \dd 2}  \big\{ \ell_s^2    (a_s^{\gg^2}-a_s^{\gg^1})_{ij} \big\} \big\|_\infty \int_{\R^d} \big|\pp_i\pp_j (1-\DD)^{-\ff \dd 2}  p_{s,t}^{\gg^1}(\cdot,y)  (x)\big|\d x\\
 &\le    c_4 \big\|\ell_s^2 (a_s^{\gg^2}-a_s^{\gg^1})_{ij}\big\|_{C_b^{\bb\land\aa}}  (t-s)^{ \ff \dd 2-1} \|\gg_s^1-\gg_s^2\|_\infty\\
 &\le c_5  \|\ell_{\nu_2}\|_{C_b^\bb}   (t-s)^{ \ff \dd 2-1} \|\gg_s^1-\gg_s^2\|_\infty. \end{align*}
By combining this with \eqref{MC4}, \eqref{MMC} and \eqref{NW*1}, we arrive at
  \beg{align*}&  \|\ell_t^1 -\ell_t^2\|_\infty \le \|\ell_{\nu_1}-\ell_{\nu_2}\|_\infty\\
  &+ c_3  \|\ell_{\nu_2}\|_\infty \bigg(\int_0^t (t-s)^{-\ff{q_0}{q_0-1} (\ff 1  2+\ff d {2p_0}) }   \|\ell_s^1-\ell_s^2\|_\infty^{\ff{q_0}{q_0-1}}\d s\bigg)^{\ff {q_0-1} {q_0 }}\\
 &+ \ff{d^2c_5} 2 \|\ell_{\nu_2}\|_{C_b^\bb}  \int_0^t    (t-s)^{\ff \dd 2 -1}
  \|\gg_s^1-\gg_s^2\|_\infty  \d s,\ \ t\in [0,T].\end{align*}
 Consequently, for any $\ll>0$,
\beg{align*} & d_{\infty,\ll}(\ell_{X^{\gg^1,\nu_1}},  \ell_{X^{\gg^2,\nu_2}}):=\sup_{t\in [0,T]} \e^{-\ll t} \|\ell_t^1-\ell_t^2\|_\infty\\
& \le \|\ell_{\nu_1}-\ell_{\nu_2}\|_\infty+\vv(\ll)  \big\{d_{\infty,\ll}(\ell_{X^{\gg^1,\nu_1}},  \ell_{X^{\gg^2,\nu_2}})+ d_{\infty,\ll}(\gg^1,\gg^2)\big\}\end{align*}
holds for
\beg{align*} \vv(\ll):= \sup_{t\in [0,T]} \bigg\{&c_3\|\ell_{\nu_2}\|_\infty\bigg(\int_0^t (t-s)^{-\ff{q_0}{q_0-1} (\ff 1  2+\ff d {2p_0}) } \e^{-\ff{q_0 \ll(t-s)}{q_0-1}}\d s\bigg)^{\ff {q_0-1}{q_0} }\\
&\qquad \qquad+ \ff{d^2c_5} 2 \|\ell_{\nu_2}\|_{C_b^\bb} \int_0^t (t-s)^{ \ff \dd 2-1} \e^{-\ll(t-s)}\d s\bigg\}.\end{align*}
Since $(p_0,q_0)\in\scr K$ implies $ \ff{q_0}{q_0-1} (\ff 1  2+\ff d {2p_0})<1$, and since $1-\ff\dd 2<1$,  by taking large enough $\ll>0$ increasing in $\|\ell_{\nu_2}\|_{C_b^\bb}$, we obtain
  \beq\label{MC6}  d_{\infty,\ll}(\ell_{X^{\gg^1,\nu_1}},  \ell_{X^{\gg^2,\nu_2}})
  \le \|\ell_{\nu_1}-\ell_{\nu_2}\|_\infty+\ff 1 4 \big\{d_{\infty,\ll}(\ell_{X^{\gg^1,\nu_1}},  \ell_{X^{\gg^2,\nu_2}})+ d_{\infty,\ll}(\gg^1,\gg^2)\big\}.\end{equation}
 Taking $\nu_1=\nu_2=\nu$, we see that the map $\gg\mapsto \ell_{X^{\gg,\nu}}$ is contractive on the complete metric space  $(L_\infty^\infty\cap \D_+^1, d_{\infty,\ll})$, so that it has a unique fixed point.
 Therefore, \eqref{E1} is well-posed. Estimate \eqref{A1} follows from Lemma \ref{L2} for $\gg_t=\ell_{X_t}$ for the solution to \eqref{E1}, while \eqref{A2} follows from
 \eqref{MC6} for $\gg_t^i:=\ell_{X_t^i}, \nu_i=\L_{X_0^i}, i=1,2.$

 \end{proof}

 \section{Density dependent reflecting SDEs}

 In this section, we extend Theorem \ref{T1}   to  density dependent reflecting SDEs on a domain $D$.  There exists additional difficulty to extend Theorem \ref{T2}, for instance,
  in the proof of Theorem \ref{T2} we used
  $$(1-\DD)^{-\ff \dd 2}\pp_i\pp_j = \pp_i\pp_j (1-\DD)^{-\ff \dd 2}$$ which is no longer true for the Nuemann Laplacian in a domain.

 Let $D\subset \R^d$ be a connected $C^2$-smooth open domain.  Consider the following density dependent reflecting SDE on the closure $\bar D$ of $D$:
\beq\label{E1'} \d X_t= b_t(X_t,\ell_{X_t}(X_t),\ell_{X_t}) \d t+ \si_t(X_t )\d W_t + \n(X_t) \d l_t,\ \ t\in [0,T],\end{equation}
where     $\n$ is the  unit inward normal vector field on the boundary $\pp D$, $l_t$ is a continuous adapted increasing process with $\d l_t$ supported on $\{t: X_t\in\pp D\},$  and
$$b: [0,T]\times \bar D\times [0,\infty)\times \D_+^1 \to \R^d,\ \ \si: [0,T]\times\bar D  \to \R^d\otimes\R^m$$ are measurable.
 Here and in the following, $\D_+^1, \tt L_q^p, L_q^p, L^p, \tt L^p,   \tt{\scr P}^p, \scr P^p$ and $\scr P$ are defined  as before for $\bar D$ replacing $\R^d$.

 We will assume $\pp D\in C_b^{2,L}$ in the following sense:
there exists a constant $r_0>0$ such that  the polar coordinate map
$$ \Psi: \pp D\times [-r_0,r_0] \ni (\theta, r)\mapsto   \theta+ r\n(\theta) \in \pp_{\pm r_0}D:=\big\{x\in\R^d: \rr_\pp(x):={\rm dist}(x,\pp D)\le r_0\big\} $$
is a $C^2$-diffeomorphism, such that   $ \Psi^{-1}(x)$  have   bounded and continuous first and second order derivatives in $x\in  \pp_{\pm r_0}D$,   and  $\nn^2\rr_\pp$ is Lipschitz continuous on $\pp_{\pm r_0}D.$

Note that $\pp D\in C_b^{2,L}$ does not imply the boundedness of $D$ or $\pp D$, but any bounded $C^{2,L}$ domain  satisfies $\pp D\in C_b^{2,L}$.

\beg{defn} (1) A pair $(X_t, l_t)_{t\in [0,T]}$ is called a  (strong) solution of \eqref{E1'},  if  $(X_t)_{t\in [0,T]}$ is a continuous adapted process on $\bar D$,
$(l_t)_{t\in [0,T]}$ is a continuous adapted increasing process with $l_0=0$ and $\d l_t$ supported on $\{t\in [0,T]: X_t\in\pp D\}$,  such that
$$ \int_0^T \E\big[|b_s(X_s,\ell_{X_s}(X_s), \ell_{X_t})|+\|\si_s(X_s)\|^2\big]  \d s<\infty $$ and
$\P$-a.s.
$$X_t=X_0+ \int_0^t b_s(X_s,\ell_{X_s}(X_s), \ell_{X_s})\d s +\int_0^t \si_s(X_s)\d W_s+\int_0^t\n(X_s)\d l_s, \ \  t\in [0,T].$$

(2) A triple  $(X_t,l_t,W_t)_{t\in [0,T]}$ is called a weak solution of \eqref{E1'}, if $(W_t)_{t\in [0,T]}$ an $m$-dimensional Brownian under a complete filtration probability space
$(\OO, \{\F_t\}_{t\in [0,T]}, \P)$ such that $(X_t,l_t)_{t\in [0,T]}$ solves $\eqref{E1'}$. We identify   any two weak solutions $(X_t, l_t,W_t)$ and $(\bar X_t,\bar l_t,\bar W_t)$ if
$(X_t, l_t)_{t\in [0,T]}$ and $(\bar X_t,\bar l_t)_{t\in [0,T]}$  have the same   distribution under  the corresponding probability spaces.
\end{defn}

To extend assumption {\bf (A)} to the present setting, we introduce  the  the Neumann semigroup  $\{P_{s,t}^{a,b^{(1)}}\}_{0\le s\le t\le T}$    generated by $  L_t^{a, b^{(1)}}$ on $\bar D$
for $a_t:=\si_t\si_t^*$,
   that is, for any $\phi\in C_b^2(\bar D)$,  and any  $t\in (0,T]$, $(P_{s,t}^{a,b^{(1)}}\phi)_{s\in [0,t]}$ is the unique solution of the PDE
\beq\label{NMM} \pp_s u_s= - L_s^{a,b^{(1)}}  u_s,\ \ \nn_{\n}u_s|_{\pp D}=0\ \text{for}\ s\in [0,t), u_t=\phi.\end{equation}
For any $t>0,$ let $C_b^{1,2}([0,t]\times\bar D)$ be the set of  functions  $f\in C_b([0,t]\times\bar D)$  with bounded and continuous derivatives $\pp_t f, \nn f$ and $\nn^2 f$.

We now extend {\bf (A)} to the domain setting as follows.

 \beg{enumerate} \item[{\bf (C)}] Let $k\in [1,\infty]$, $\pp D\in C_b^{2,L}$,   $\si_t(x,\rr)=\si_t(x)$ and $b_t(x,r,\rr)= b_t^{(1)}(x)+ b_t^{(0)}(x,r,\rr)$ satisfy the following conditions.
 \item [$(C_1)$] $a_t(x):=(\si_t\si_t^*)(x)$ is invertible for $(t,x)\in [0,T]\times\bar D$, $\|a\|_\infty+\|a^{-1}\|_\infty<\infty$,  and
 $$\lim_{\vv\downarrow 0} \sup_{t\in [0,T]}\sup_{x,y\in \bar D,|x-y|\le \vv} \|a_t(x)-a_t(y)\|=0.$$
 Moreover, $\si_t$ is weakly differentiable with $\|\nn\si\|\le \sum_{i=1}^l f_i$ for some
 $l\in \mathbb N$,  $0\le f_i\in \tt L_{q_i}^{p_i}, (p_i,q_i)\in \scr K, 1\le i\le l$.
 \item[$(C_2)$]  $(A_2)$ holds for $\bar D$ replacing $\R^d$.
 \item[$(C_3)$]  For any  $\phi\in C_b^2(\bar D)$ and $t\in (0,T]$, the PDE \eqref{NMM} has a unique solution $P_{\cdot,t}^{\si,b^{(1)}}\phi\in C_b^{1,2}([0,t]\times \bar D),$  such that    for some  constants $c,\kk>0$ and diffeomorphsims $\{\psi_{s,t}\}_{0\le s\le t\le T}$ on $\bar D$ satisfying $\eqref{MP1}$, the hear kernel $p_{s,t}^{a, b^{(1)}}$ of $P_{s,t}^{a,b^{(1)}}$ satisfies
 \beq\label{AB1}  |\nn^i p_{s,t}^{a,b^{(1)}}(\cdot,y)|(x)  \le  c (t-s)^{-\ff i 2} p_{s,t}^\kk(\psi_{s,t}(x)-y),\ \ 0\le s<t\le T, x,y\in \bar D,\ i=1,2.\end{equation}
 \end{enumerate}

 By  \cite[Theorem VI.3.1]{Carr}, $(C_3)$ holds   if  $D$ is  bounded with $\pp D\in C^{2+\aa}$ for some $\aa\in (0,1),$  and    there exists   $c>0$ such that
$$ \big\{|b_t^{(1)}(x)-b_s^{(1)}(y)|+\|a_t(x)-a_s(y)\|\big\}\le c (|t-s|^{\aa}+|x-y|^{\aa}),\ \ s,t\in [0,T],x,y\in \bar D.$$
If moreover $\nn a_t$ is H\"older continuous uniformly in $t\in [0,T]$, then for any $\bb\in (0,1)$ there exists a constant $c>0$ such that \eqref{MP3''} holds for $p_{s,t}^{a, b^{(1)}}$:
\beq\label{HD}  \beg{split} & |\nn p_{s,t}^{a,b^{(1)}}(\cdot,y)(x)-\nn p_{s,t}^{a,b^{(1)}}(\cdot,y')(x) |\\
&\le c |y-y'|^\bb (t-s)^{-\ff {1+\bb} 2} \big\{p_{t-s}^\kk(\psi_{s,t}(x)-y)+ p_{t-s}^\kk(\psi_{s,t}(x)-y')\big\},\\
&\qquad  \  \ 0\le s<t\le T, \ x,x',y\in\R^d.\end{split}\end{equation}
The following result extends Theorem  \ref{T1} to the reflecting setting.

\beg{thm}\label{T3} Assume {\bf (C)}   for some  $k\in [\ff{p_0}{p_0-1},\infty]\cap (k_0,\infty],$ where $k_0:=\ff d {2\theta+1-d p_0^{-1}-2 q_0^{-1}}$.
 \beg{enumerate}
\item[$(1)$] For any   $\nu\in \tt {\scr P}^k,$ $\eqref{E1'}$ has a unique strong $($respectively   weak$)$ solution with $\L_{X_0}=\nu$ satisfying $\ell_{X_\cdot}\in \tt L_\infty^k(\bar D)$.
 Moreover, there exists a constant $c>0$ such that for  any two   solutions $X_t^1$ and $ X_t^2$ of $\eqref{E1'}$ with   initial
 distributions  $\L_{X_\cdot^1}, \L_{X_\cdot^2}\in\tt L^k_\infty$,
 $$ \sup_{t\in [0,T]} \|\ell_{X_t^1} -\ell_{X_t^2}\|_{\tt L^k}\le c \|\L_{X_0^1}- \L_{X_0^2}\|_{\tt L^k}.$$
 \item[$(2)$] Assertions in $(1)$ hold  for $(\scr P^k, L_\infty^k, L^k)$ replacing $(\tt{\scr P}^k, \tt L_\infty^k, \tt L^k)$, provided in $(C_2)$  the condition $(A_2)$ is replaced by  $(A_2')$ for $\bar D$ replacing $\R^d$. \end{enumerate}
  \end{thm}

   \beg{proof}   As explained  in the proof of Theorem \ref{T1}(2), we only    prove the first assertion.

   According to \cite[Theorem 2.2(ii)]{W21b}, for any $\gg\in \tt L_\infty^k \cap \D_+^1$, the reflecting SDE
\beq\label{E2'} \d  X_t^\gg = b_t^\gg(X_t^\gg)\d t+\si_t(X_t^\gg)\d W_t + \n(X_t^\gg)\d l_t^\gg,\ \ t\in [0,T]    \end{equation}
is well-posed. Let $X_t^{\gg,x}$ denote the solution with initial value $X_0^\gg=x\in \bar D$, and simply denote $X_t^\gg$ for the solution with $X_0^\gg=X_0$ for $\L_{X_0}=\nu.$

By Theorem 6.2.7(ii)-(iii) in \cite{BKRS}, the distribution density function $\ell_{X_t^{\gg,x}}$ exists for $t\in (0,T]$ and $x\in D$.
  Next, by \cite[Theorem 4.1]{W21b} for distribution independent drift, there exists $c>0$ such that the following log-Harnack inequality holds for
  the associated semigroup:
  $$P_t^\gg \log f(x)\le \log P_t f(y)+\ff{c|x-y|^2}{t},\ \ t\in (0,T], x,y\in \R^d, f>0.$$
  This implies that $\{\L_{X_t^{\gg,x}}\}_{x\in\bar D}$ are mutually equivalent for $t\in (0,T]$. Thus,  the existence of $\{\ell_{X_t^{\gg,x}}\}_{t\in (0,T]}$ for $x\in D$ implies that  for $x\in \bar D$. Consequently,
  $$\Phi_t^\nu \gg:=\ell_{X_t^\gg}=\int_{\bar D} \ell_{X_t^{\gg,x}}\L_{X_0}(\d x)$$ exists for any $t\in (0,T].$
Since $\{\psi_{s,t}\}_{0\le s\le t\le T}$ are diffeomorphisms on $\bar D$ satisfying \eqref{MP1},
$$\hat P_{s,t}^\kk f(y):=\int_{\bar D} p_{t-s}^\kk(y-\psi_{s,t}(x)) f(x)\d x,\ \ y\in \bar D$$
gives rise to a family of linear operators satisfying \eqref{M0} and Lemma \ref{LA} for norms defined with $\bar D$ replacing $\R^d$. Then by repeating the proof  of Lemma \ref{L1}  using the present estimates, we conclude that $\Phi^\nu$ maps $\tt{\scr P}^k_{\nu,T}$ into  $\tt{\scr P}^k_{\nu,T}$ such that \eqref{ETT} and \eqref{*2} hold  for $ \bar D $ replacing $ \R^d.$ With this result and using $(C_3)$ replacing \eqref{MP3}, we  prove Theorem \ref{T3}(1) by the  means in
used the proof of Theorem \ref{T1}(1).
\end{proof}

  \paragraph{Acknowledgement.} The author would like to thank the referee for helpful comments and suggestions.

\end{document}